\documentclass[11pt]{amsart}

\usepackage{amsmath,amssymb,amsthm,mathtools}
\usepackage{enumitem}
\usepackage{microtype}
\usepackage[colorlinks=true,linkcolor=blue,citecolor=blue,urlcolor=blue]{hyperref}
\usepackage[margin=1.15in]{geometry}

\numberwithin{equation}{section}
\newtheorem{theorem}{Theorem}[section]
\newtheorem{lemma}[theorem]{Lemma}
\newtheorem{corollary}[theorem]{Corollary}
\newtheorem{proposition}[theorem]{Proposition}
\newtheorem*{fernandez}{Fern\'andez's theorem}
\theoremstyle{definition}
\newtheorem{example}[theorem]{Example}
\newtheoremstyle{remarkstyle}
  {}{}{\normalfont}{}{\itshape}{.}{ }{}
\theoremstyle{remarkstyle}
\newtheorem*{remark}{Remark}
\allowdisplaybreaks[4]
\newcommand{\D}{\mathbb{D}}
\newcommand{\T}{\mathbb{T}}
\newcommand{\R}{\mathbb{R}}
\newcommand{\C}{\mathbb{C}}
\newcommand{\Cbar}{\overline{\mathbb{C}}}
\newcommand{\Rbar}{\overline{\mathbb{R}}}
\newcommand{\Hp}{\mathbb{H}^{+}}
\newcommand{\Hm}{\mathbb{H}^{-}}
\newcommand{\BMOA}{\mathrm{BMOA}}
\newcommand{\BMO}{\mathrm{BMO}}
\newcommand{\VMO}{\mathrm{VMO}}
\newcommand{\VMOA}{\mathrm{VMOA}}
\newcommand{\CM}{\mathrm{CM}}
\newcommand{\Aut}{\operatorname{Aut}}
\newcommand{\diam}{\operatorname{diam}}
\newcommand{\dist}{\operatorname{dist}}
\newcommand{\length}{\operatorname{length}}
\newcommand{\artanh}{\operatorname{artanh}}
\renewcommand{\Re}{\operatorname{Re}}
\renewcommand{\Im}{\operatorname{Im}}

\hypersetup{
  pdftitle={$\mathrm{BMOA}$ and Carleson-homogeneous Denjoy domains},
  pdfauthor={Katsuhiko Matsuzaki and Fei Tao}
}

\title[$\mathrm{BMOA}$ and Carleson-homogeneous Denjoy domains]
{$\mathrm{BMOA}$ and Carleson-homogeneous Denjoy domains}

\author{Katsuhiko Matsuzaki}
\address{Department of Mathematics, School of Education, Waseda University, Tokyo 169-8050, Japan}
\email{matsuzak@waseda.jp}
\author{Fei Tao}
\address{Beijing International Center for Mathematical Research, Peking University, Beijing 100871, P. R. China}
\email{ferrytau@pku.edu.cn}


\thanks{The first author is partially supported by Japan Society for the Promotion of Science (KAKENHI 23K25775 and 23K17656); the second author is partially supported by National Key R \& D Program of China (2025YFA1017500)}
\date{}

\begin{document}

\begin{abstract}
Let $E$ be a proper closed subset of the extended real line that contains $[-\infty,-1]\cup[1,+\infty]$, and let $F$ be a universal covering of the Denjoy domain $\Omega=\Cbar\setminus E$. We prove the conjecture of Huo and Zinsmeister that $E$ is Carleson-homogeneous if $\log F'\in\BMOA(\D)$. The proof is based on an estimate for the length of $F^{-1}(\R)$ in terms of the Bloch seminorm of $\log F'$ and a weighted $L^2$ integral of $F''/F'$. From this estimate we deduce that $\Omega$ has the Hayman--Wu property, and Carleson-homogeneity then follows from a theorem of Fern\'andez. We also characterize these domains by Carleson measure conditions for the arclength on $F^{-1}(\R)$ and for the Schwarzian derivative of $F$, and by the $H^\infty$ interpolation property of the inverse images of suitable sequences in $\R\setminus E$.
\end{abstract}

\maketitle

\section{Introduction}\label{sec:intro}

A Denjoy domain is a plane domain whose complement is contained in the real line $\mathbb{R}$. Throughout this paper, $E$ denotes a closed subset of the extended real line $\Rbar=\R\cup\{\infty\}$ such that
\begin{equation}\label{eq:exterior}
[-\infty,-1]\cup[1,+\infty]\subset E,\qquad \Rbar\setminus E\neq\emptyset ,
\end{equation}
where $-\infty$ and $+\infty$ denote the same point $\infty$. Set $\Omega=\Cbar\setminus E$. Since $\R\setminus E\ne\emptyset$, the upper and lower half-planes are joined through the components of $\R\setminus E$, so $\Omega$ is a connected hyperbolic Denjoy domain with $\infty\in\partial\Omega$. Let $F\colon\D\to\Omega$ be a universal covering and let $G$ be its Fuchsian group, so that $F\circ\gamma=F$ for $\gamma\in G$ and $\Omega\cong\D/G$. The derivative $F'$ is zero-free, and $\log F'$ denotes an analytic branch of its logarithm. All conditions on $\log F'$ below are independent of the choice of the branch. Any two universal coverings of $\Omega$ differ by precomposition with an automorphism of $\D$.

We say that $E$ is \emph{Carleson-homogeneous}, and that $\Omega$ is a \emph{Carleson-homogeneous Denjoy domain}, if there is a constant $C>0$ such that
\begin{equation}\label{eq:homogeneity}
|E\cap(x-r,x+r)|\ge C\,r\qquad(x\in E\cap\R,\ r>0),
\end{equation}
where $|\cdot|$ denotes Lebesgue measure on $\R$. Such sets were introduced by Carleson \cite{Ca} in his work on $H^\infty$ in complements of subsets of the line,
and have subsequently played an important role in the study of Hardy spaces \cite{Z}, harmonic measure, and interpolation \cite{JM}.
%and they appear in several characterizations of Denjoy domains, for instance in terms of Hardy spaces \cite{Z} and of harmonic measure \cite{JM}.

Let $\mathcal B(\D)$ be the Bloch space, with seminorm $\|\varphi\|_{\mathcal B}=\sup_{z\in\D}(1-|z|^2)|\varphi'(z)|$. The space $\BMOA(\D)\subset\mathcal B(\D)$ is recalled in Section~\ref{sec:prelim}. For simply connected domains, the Riemann maps $f$ with $\log f'\in\BMOA(\D)$ were characterized by Bishop and Jones \cite{BJ}. Huo and Zinsmeister \cite[Theorem~1.1]{HZ} proved that $\log F'\in\BMOA(\D)$ for Carleson-homogeneous Denjoy domains and conjectured the converse. We prove this converse.

For a set $X\subset\C$, $\length(X)$ denotes its one-dimensional Euclidean Hausdorff measure. For a rectifiable Jordan arc or Jordan curve, this agrees with its usual arclength. By arclength on $X$ we mean the restriction of this measure to $X$, denoted by $ds|_X$.

The Hayman--Wu theorem gives a universal bound for the length of inverse images of lines under conformal maps from $\D$ onto simply connected proper planar domains \cite{HW}. The proof below uses the corresponding property for universal coverings.

We say that a hyperbolic domain $V\subset\Cbar$ has the \emph{Hayman--Wu property} if there is a constant $C_V>0$ such that
\begin{equation}\label{eq:HW-intro}
\length\bigl(F^{-1}(L\cap V)\bigr)\le C_V
\end{equation}
for every universal covering $F$ of $V$ and every line $L\subset\C$. When $\infty\in V$, the universal covering $F$ is meromorphic, and its poles do not belong to $F^{-1}(L\cap V)$. For a compact set $K\subset\R$ of positive diameter, Carleson-homogeneity is understood with the obvious restriction on centers and scales,
\begin{equation}\label{eq:compact-homogeneity}
|K\cap(y-s,y+s)|\ge C\,s\qquad(y\in K,\ 0<s\le\diam K),
\end{equation}
and this condition is invariant under real affine maps.

\begin{fernandez}[{\cite[Theorem~4]{Fe}}]
Let $K\subset\R$ be a compact set such that $V=\Cbar\setminus K$ is hyperbolic. Then $V$ has the Hayman--Wu property if and only if $K$ is Carleson-homogeneous in the sense of \eqref{eq:compact-homogeneity}.
\end{fernandez}

The passage from compact sets to the normalization \eqref{eq:exterior} is made in Section~\ref{subsec:HW}.

The Schwarzian derivative of $F$ is $S_F=(F''/F')'-\frac12(F''/F')^2$. Carleson measures on $\D$, whose class is denoted by $\CM(\D)$, are recalled in Section~\ref{sec:prelim}.

\begin{theorem}\label{thm:main}
Let $E$ satisfy \eqref{eq:exterior}, and let $F$ be a universal covering of $\Omega=\Cbar\setminus E$. The following conditions are equivalent.
\begin{enumerate}[label=\textup{(\roman*)}]
\item $\Omega$ is a Carleson-homogeneous Denjoy domain.
\item $\log F'\in\BMOA(\D)$.
\item Arclength on $F^{-1}(\R)$ is a Carleson measure on $\D$.
\item $\Omega$ has the Hayman--Wu property.
\item $(1-|z|^2)^3|S_F(z)|^2\,dxdy\in\CM(\D)$.
\end{enumerate}
\end{theorem}

Conditions (ii), (iii) and (v) do not depend on the choice of $F$. The implication (i)$\Rightarrow$(ii) is \cite[Theorem~1.1]{HZ}, and the equivalence (iv)$\Leftrightarrow$(i) follows from Fern\'andez's theorem after the normalization in Section~\ref{subsec:HW}. The equivalence (ii)$\Leftrightarrow$(v) is known (see \cite{AZ,GNR}), and Section~\ref{sec:schwarzian} gives a direct proof for locally univalent analytic functions.

The main new step is (ii)$\Rightarrow$(iii). The implication (iii)$\Rightarrow$(iv) follows the real-line reduction of Fern\'andez \cite[Section~4, p.~54]{Fe}.
Section~\ref{sec:circles} gives the details and extends the estimates to circles. 

We first prove the following length estimate.

\begin{theorem}\label{thm:length}
Let $E\subset\Rbar$ be closed and satisfy \eqref{eq:exterior}, and let
$F\colon\D\to\Omega=\Cbar\setminus E$ be a universal covering. Then
\begin{equation}\label{eq:length-estimate}
\length\bigl(F^{-1}(\R)\bigr)\le C\bigl(1+\|\log F'\|_{\mathcal B}\bigr)^2\Bigl(1+\int_{\D}\Bigl|\frac{F''(z)}{F'(z)}\Bigr|^2(1-|z|^2)\,dxdy\Bigr),
\end{equation}
where $C>0$ is an absolute constant and the right-hand side is allowed to be infinite.
\end{theorem}

The integral in \eqref{eq:length-estimate} is finite if and only if $\log F'\in H^2(\D)$ (Lemma~\ref{lem:bmo-basic}). Thus $F^{-1}(\R)$ has finite length whenever $\log F'\in H^2(\D)$ and $\Omega$ is uniformly perfect in the sense recalled in Section~\ref{subsec:hyperbolic}. Every component of $F^{-1}(\R)$ is a complete hyperbolic geodesic, and a bound for $\|\log F'\|_{\mathcal B}$ forces distinct components to be uniformly separated in the hyperbolic metric. Since $F$ maps a circular component of Euclidean radius $r$ into a straight line, the curvature of this component gives $|F''/F'|\ge1/r$ on it. A lower bound of this type is suggested in the final remark of \cite{HZ}. Theorem~\ref{thm:length} then follows by integrating $|F''/F'|^2$ over disjoint disks centered on the components.

If $\log F'\in\BMOA(\D)$, the Garsia seminorm of $\log(F\circ T)'$ is bounded independently of the automorphism $T$ of $\D$, and Theorem~\ref{thm:length} applies uniformly to the universal coverings $F\circ T$. The components of $\D\setminus F^{-1}(\R)$ are mapped conformally onto half-planes, and their boundaries are rectifiable Jordan curves. The Carleson embedding theorem for $H^1(\D)$ then gives a bound for the inverse image of every circle and every line.

\begin{theorem}\label{thm:uniform}
Suppose that $\log F'\in\BMOA(\D)$, and let $\|\log F'\|_*$ be its Garsia seminorm, defined in \eqref{eq:garsia}. Then
\begin{equation}\label{eq:uniform-circles}
\length\bigl((F\circ T)^{-1}(\mathcal C)\bigr)\le C\bigl(1+\|\log F'\|_*\bigr)^4
\end{equation}
for every automorphism $T$ of $\D$ and every circle or line $\mathcal C$ in $\C$, where $C$ is an absolute constant. Moreover, arclength on $F^{-1}(\R)$ is a Carleson measure on $\D$ with Carleson norm at most $C(1+\|\log F'\|_*)^4$.
\end{theorem}

Under the hypothesis of Theorem~\ref{thm:uniform}, choose $a\in\R\setminus E$ and set $\Phi(w)=1/(w-a)$, so that $K=\Phi(E)$ is a compact subset of $\R$. The circle and line estimates imply that $\Cbar\setminus K$ has the Hayman--Wu property. Fern\'andez's theorem and Lemma~\ref{lem:transfer} then show that $K$, and hence $E$, is Carleson-homogeneous. This argument does not provide an explicit lower bound for the homogeneity constant in terms of $\|\log F'\|_*$.

We prove two supplemental results around Theorem~\ref{thm:main}. The first concerns the little Bloch condition.
Although $\log F'$ belongs to the Bloch space $\mathcal B(\D)$ if and only if $\Omega$ is uniformly perfect
by Pommerenke \cite{Po84},
we can prove that $(1-|z|^2)|F''(z)/F'(z)|$ does not tend to zero as $|z|\to1$. In particular, $\log F'\notin\VMOA(\D)$ for every universal covering of every domain $\Omega$ as above (Proposition~\ref{prop:little-bloch}).
The other investigates interpolating sequences on $\D$ for $H^\infty$.
For sequences $Z\subset\R\setminus E$ satisfying the separation and covering conditions in Theorem~\ref{thm:interp-char}, we prove that the set $E$ is Carleson-homogeneous 
if and only if $F^{-1}(Z)$ is an interpolating sequence for $H^\infty$. 

The organization of this paper is as follows.
Section~\ref{sec:prelim} collects the preliminaries, and Section~\ref{sec:length} proves the length estimate and its uniform version.
Section~\ref{sec:circles} proves the equivalence of (i)--(iv).
Section~\ref{sec:schwarzian} completes the proof of Theorem~\ref{thm:main} through the Schwarzian condition and concludes with the little Bloch space. Finally, in Section~\ref{sec:interpolation}, we relate interpolating sequences on $\D$ for $H^\infty$ to the equivalent conditions in Theorem~\ref{thm:main}.

Throughout, $C$ denotes a positive absolute constant whose value may change from line to line, unless a dependence on parameters is indicated. Numbered constants such as $C_1$ and $C_S$ are fixed.

\section{Preliminaries}\label{sec:prelim}

\subsection{Hardy spaces and BMOA}

Let $\T=\partial\D$, and let $dm(\zeta)=|d\zeta|/(2\pi)$ be normalized arclength on $\T$. For a Borel set $X\subset\T$ we also write $|X|$ for its unnormalized arclength, so that $|X|=2\pi m(X)$. For $0<p<\infty$ and $u$ analytic in $\D$ we write
\[
\|u\|_{H^p}=\sup_{0<s<1}\Bigl(\int_{\T}|u(s\zeta)|^p\,dm(\zeta)\Bigr)^{1/p},
\]
$H^p(\D)$ is the space of analytic functions $u$ in $\D$ with $\|u\|_{H^p}<\infty$, and $u^*$ denotes the almost everywhere non-tangential boundary values of $u\in H^p(\D)$. For any $a\in\D$, we define $\sigma_a(z)=(a-z)/(1-\overline az)$. This conformal automorphism of $\D$ interchanges $0$ and $a$ and satisfies $\sigma_a\circ\sigma_a=\mathrm{id}$. Every automorphism $T$ of $\D$ can be written as $T=\sigma_{T(0)}\circ R$ with a rotation $R$, and the group of automorphisms is denoted by $\Aut(\D)$.

For $v\in L^1(\T)$ and an arc $I\subset\T$ of positive length, let $v_I=m(I)^{-1}\int_I v\,dm$. We say that $v\in\BMO(\T)$ if
\[
\|v\|_{\BMO}
=\sup_{I\subset\T}\frac{1}{m(I)}\int_I |v-v_I|\,dm<\infty,
\]
where the supremum is taken over all arcs $I\subset\T$. We say that $v\in\VMO(\T)$ if $v\in\BMO(\T)$ and
\[
\frac{1}{m(I)}\int_I |v-v_I|\,dm\longrightarrow0
\qquad\text{as }|I|\longrightarrow0.
\]

Recall that an analytic function $\varphi$ in $\D$ belongs to $\BMOA(\D)$ if $\varphi\in H^2(\D)$ and its boundary function $\varphi^*$ belongs to $\BMO(\T)$. Similarly, $\varphi\in\VMOA(\D)$ if $\varphi\in H^2(\D)$ and $\varphi^*\in\VMO(\T)$.

For an analytic function $\varphi$ in $\D$, set
\begin{equation}\label{eq:garsia}
\|\varphi\|_*=
\sup_{a\in\D}
\|\varphi\circ\sigma_a-\varphi(a)\|_{H^2}.
\end{equation}
The quantity $\|\varphi\|_*$ is finite if and only if $\varphi\in\BMOA(\D)$. On $\BMOA(\D)$ it is the Garsia seminorm, comparable to the $\BMO$ seminorm of $\varphi^*$ \cite[Chapter~VI]{Ga}. For $\varphi\in H^2(\D)$, a change of variables using $|\sigma_a'(\zeta)|=(1-|a|^2)/|\zeta-a|^2$ on $\T$ gives
\[
\|\varphi\circ\sigma_a-\varphi(a)\|_{H^2}^2=\int_{\T}|\varphi^*(\zeta)-\varphi(a)|^2\,\frac{1-|a|^2}{|\zeta-a|^2}\,dm(\zeta).
\]

\begin{lemma}\label{lem:bmo-basic}
Let $\varphi$ be analytic in $\D$.
\begin{enumerate}[label=\textup{(\alph*)}]
\item If $\varphi\in\BMOA(\D)$ and $T\in\Aut(\D)$, then $\|\varphi\circ T\|_*=\|\varphi\|_*$.
\item If $\varphi\in\BMOA(\D)$, then $\|\varphi\|_{\mathcal B}\le\|\varphi\|_*$.
\item If $\varphi(z)=\sum_{n\ge0}c_nz^n$, then
\begin{equation}\label{eq:area-identity}
\int_{\D}|\varphi'(z)|^2(1-|z|^2)\,dxdy=\pi\sum_{n\ge1}\frac{n}{n+1}|c_n|^2 .
\end{equation}
Consequently
\[
\frac{\pi}{2}\|\varphi-\varphi(0)\|_{H^2}^2\le\int_{\D}|\varphi'(z)|^2(1-|z|^2)\,dxdy\le\pi\|\varphi-\varphi(0)\|_{H^2}^2.
\]
When $\varphi\in\BMOA(\D)$, the right-hand side is at most $\pi\|\varphi\|_*^2$.
\end{enumerate}
\end{lemma}

\begin{proof}
(a) For $a\in\D$, the automorphisms $T\circ\sigma_a$ and $\sigma_{T(a)}$ both map $0$ to $T(a)$, so that $T\circ\sigma_a=\sigma_{T(a)}\circ R_a$ for a rotation $R_a$. The $H^2$ norm is invariant under rotations, and hence
\[
\|\varphi\circ T\circ\sigma_a-\varphi(T(a))\|_{H^2}=\|\varphi\circ\sigma_{T(a)}-\varphi(T(a))\|_{H^2}.
\]
Taking the supremum over $a\in\D$ proves (a).

(b) The coefficient of $z$ in the Taylor expansion of $\varphi\circ\sigma_a-\varphi(a)$ is $-(1-|a|^2)\varphi'(a)$, and the modulus of each Taylor coefficient is bounded by the $H^2$ norm.

(c) By orthogonality on circles,
\[
\int_{\D}|\varphi'(z)|^2(1-|z|^2)\,dxdy=2\pi\sum_{n\ge1}n^2|c_n|^2\int_0^1s^{2n-1}(1-s^2)\,ds=\pi\sum_{n\ge1}\frac{n}{n+1}|c_n|^2 .
\]
The two inequalities follow from $\frac12\le n/(n+1)<1$, and the last assertion is the case $a=0$ of \eqref{eq:garsia}, since $\sigma_0(z)=-z$.
\end{proof}

\begin{lemma}\label{lem:log-aut}
Let $T\in\Aut(\D)$ and $k_T=\log T'-\log T'(0)$. Then $\|k_T\|_{H^2}^2\le2\pi^2/3$ and $\|\log T'\|_*\le4\pi/\sqrt6$.
\end{lemma}

\begin{proof}
Write $T(z)=e^{i\theta}(z-\alpha)/(1-\overline\alpha z)$. Then
\[
k_T(z)=-2\log(1-\overline\alpha z)=2\sum_{n\ge1}\frac{\overline\alpha^{\,n}}{n}z^n,
\qquad
\|k_T\|_{H^2}^2=4\sum_{n\ge1}\frac{|\alpha|^{2n}}{n^2}\le\frac{2\pi^2}{3}.
\]
For $a\in\D$ the chain rule gives $k_T\circ\sigma_a-k_T(a)=k_{T\circ\sigma_a}-k_{\sigma_a}$. The bound for the $H^2$ norm, applied to $T\circ\sigma_a$ and to $\sigma_a$, shows that $\|k_T\circ\sigma_a-k_T(a)\|_{H^2}\le 2\sqrt{2\pi^2/3}=4\pi/\sqrt6$.
\end{proof}

\subsection{Carleson measures}

A positive Borel measure $\mu$ on $\D$ is called a \emph{Carleson measure}, and we write $\mu\in\CM(\D)$, if there is a constant $C>0$ such that $\mu(D(\zeta,r))\le Cr$ for all $\zeta\in\T$ and $r>0$, where $D(\zeta,r)$ is the Euclidean disk of center $\zeta$ and radius $r$. The infimum of all such constants is denoted by $\|\mu\|_{\CM}$. Up to absolute constants, this agrees with the usual definition by Carleson boxes. By Carleson's embedding theorem \cite[Theorem~9.3]{Du}, for every $0<p<\infty$ there is a constant $C_p$ such that
\begin{equation}\label{eq:carleson-embedding}
\int_{\D}|u|^p\,d\mu\le C_p\|\mu\|_{\CM}\|u\|_{H^p}^p\qquad(u\in H^p(\D)).
\end{equation}
Only the case $p=1$ will be used.

We shall also need the following conformally invariant form of the Carleson condition.

\begin{lemma}[{\cite[Chapter~VI, Lemma~3.3]{Ga}}]\label{lem:kernel}
For a positive Borel measure $\mu$ on $\D$, set
\[
S(\mu)=\sup_{a\in\D}\int_{\D}
\frac{1-|a|^2}{|1-\overline a z|^2}\,d\mu(z).
\]
Then $\mu\in\CM(\D)$ if and only if $S(\mu)<\infty$.
In this case,
\[
C_S^{-1}\|\mu\|_{\CM}
\le S(\mu)
\le C_S\|\mu\|_{\CM},
\]
where $C_S\ge1$ is an absolute constant.
\end{lemma}

Since $(1-|a|^2)/|1-\overline az|^2=|\sigma_a'(z)|$, the quantity $S(\mu)$ is equal to $\sup_{T\in\Aut(\D)}\int_{\D}|T'|\,d\mu$.

We recall a classical theorem of F. and M. Riesz. If $U\subset\C$ is a Jordan domain whose boundary is rectifiable and $g$ is a Riemann map of $\D$ onto $U$, then $g'\in H^1(\D)$, the boundary function of $g$ 
is absolutely continuous on $\T$, and
\begin{equation}\label{eq:rectifiable}
2\pi\|g'\|_{H^1}=\length(\partial U)
\end{equation}
(see \cite[Theorems~3.11 and~3.12]{Du}).

\subsection{Hyperbolic geometry}\label{subsec:hyperbolic}

The hyperbolic metric on $\D$ is $\lambda_\D(z)|dz|=2|dz|/(1-|z|^2)$, which has curvature $-1$, and $d_\D$ denotes the associated distance, so that $d_\D(0,z)=\log\bigl((1+|z|)/(1-|z|)\bigr)$. The hyperbolic disk $D_\D(a,R)$ of center $a$ and radius $R$ is the image of $D(0,\tanh(R/2))$ under $\sigma_a$. Complete geodesics are diameters and circular arcs orthogonal to $\T$, and the reflection in a geodesic is an anticonformal isometry of $\D$.

The hyperbolic density $\lambda_\Omega$ of $\Omega$ is determined by $\lambda_\Omega(F(z))|F'(z)|=\lambda_\D(z)$, so that the covering map $F:\D \to \Omega$ is a local isometry. Every nonidentity element of $G$ has no fixed point in $\D$. We will use the elementary inequality
\begin{equation}\label{eq:harnack}
\frac{1-|w|^2}{1-|u|^2}\le e^{d_\D(u,w)}\qquad(u,w\in\D).
\end{equation}
Indeed, if $v=\sigma_u(w)$ and $t=|v|$, then $w=\sigma_u(v)$ and
\[
\frac{1-|w|^2}{1-|u|^2}=\frac{1-t^2}{|1-\overline uv|^2}\le\frac{1-t^2}{(1-t)^2}=\frac{1+t}{1-t}=e^{d_\D(u,w)} .
\]

\begin{lemma}\label{lem:injectivity}
Let $f$ be a locally univalent analytic function in $\D$ with $M=\|\log f'\|_{\mathcal B}<\infty$, and let
\[
\varepsilon_0=\frac{1}{4(1+M)},\qquad \rho_0=2\artanh\frac{\varepsilon_0}{4}.
\]
Then $f$ is injective on every Euclidean disk $D(a,\varepsilon_0(1-|a|))$ and on every hyperbolic disk $D_\D(a,\rho_0)$ for $a\in\D$. Moreover $\rho_0\ge\varepsilon_0/2$.
\end{lemma}

\begin{proof}
Fix $a\in\D$, let $d=1-|a|$, and let $\varphi=\log f'$. If $z\in D(a,\varepsilon_0d)$, then every point $w$ of the segment $[a,z]$ satisfies $1-|w|^2\ge1-|w|\ge(1-\varepsilon_0)d$, and therefore
\[
|\varphi(z)-\varphi(a)|\le\frac{M|z-a|}{(1-\varepsilon_0)d}\le\frac{M\varepsilon_0}{1-\varepsilon_0}=\frac{M}{4M+3}<\frac14 .
\]
Since $f'(z)/f'(a)=\exp\bigl(\varphi(z)-\varphi(a)\bigr)$ and $|\Im(\varphi(z)-\varphi(a))|<1/4<\pi/2$, it follows that $\Re\bigl(f'(z)/f'(a)\bigr)>0$ on $D(a,\varepsilon_0d)$. For distinct points $z_1,z_2$ of this convex set $D(a,\varepsilon_0d)$,
\[
\frac{f(z_2)-f(z_1)}{(z_2-z_1)f'(a)}=\int_0^1\frac{f'(z_1+t(z_2-z_1))}{f'(a)}\,dt
\]
has positive real part, so that $f(z_1)\ne f(z_2)$.

The automorphism $\sigma_a$ maps $D(0,\varepsilon_0/4)$ onto $D_\D(a,\rho_0)$. For $|w|<\varepsilon_0/4$ we have
\[
|\sigma_a(w)-a|=\frac{(1-|a|^2)|w|}{|1-\overline aw|}\le\frac{2d|w|}{1-|w|}<\varepsilon_0d,
\]
because $\varepsilon_0\le1/4$. Thus $D_\D(a,\rho_0)\subset D(a,\varepsilon_0d)$. Finally, $\artanh t\ge t$ gives $\rho_0\ge\varepsilon_0/2$.
\end{proof}

In the normalization \eqref{eq:exterior}, the set $E$, and the domain $\Omega$, are called \emph{uniformly perfect} if there is $c\in(0,1)$ such that $E\cap\{w:cr\le|w-a|\le r\}\ne\emptyset$ for every $a\in E\cap\R$ and $r>0$. By Pommerenke \cite{Po84}, this is equivalent to $\log F'\in\mathcal B(\D)$, and by Beardon and Pommerenke \cite{BP} it is also equivalent to $\inf_{w\in\Omega}\lambda_\Omega(w)\dist(w,E)>0$.

\subsection{The inverse image of the real line}

Since $F$ takes values in $\Omega$, we have $F^{-1}(\R)=F^{-1}(\R\setminus E)$. The following lemma is a form of the Schwarz reflection principle.

\begin{lemma}\label{lem:reflection}
The following statements hold:
\begin{enumerate}[label=\textup{(\alph*)}]
\item Every component $\ell$ of $F^{-1}(\R)$ is a complete hyperbolic geodesic, and the reflection $R_\ell$ in $\ell$ satisfies $F\circ R_\ell=\overline F$.
\item The components of $F^{-1}(\R)$ form a countable family that is locally finite in $\D$.
\item Each component $\ell$ of $F^{-1}(\R)$ is mapped by $F$ homeomorphically onto a component of $\R\setminus E$.
\end{enumerate}
\end{lemma}

\begin{proof}
(a) The domain $\Omega$ is symmetric with respect to $\R$, so $\overline F$ takes values in $\Omega$. Let $a\in F^{-1}(\R)$. Since $F(a)$ is real, the lifting property of the universal covering gives a unique continuous map $\iota\colon\D\to\D$ with $F\circ \iota=\overline F$ and $\iota(a)=a$. Locally, $\iota$ is the composition of $\overline F$ with a branch of $F^{-1}$, so $\iota$ is antiholomorphic. Moreover $F\circ \iota^2=F$ and $\iota^2(a)=a$, and the uniqueness of lifts gives $\iota^2=\mathrm{id}$. Thus $\iota$ is an antiholomorphic automorphism of $\D$ with a fixed point. After conjugation by an automorphism that maps $a$ to $0$, it takes the form $z\mapsto e^{i\theta}\overline z$, whose fixed set is a diameter. Hence $\iota$ is the reflection in a complete geodesic $\ell_a$ through $a$, and $F=\overline F$ on $\ell_a$, that is, $\ell_a\subset F^{-1}(\R)$.

Near each of its points, $F^{-1}(\R)$ is a single smooth arc, since it is the inverse image of a line under a locally conformal map. Near a point of $\ell_a$ this arc coincides with $\ell_a$. Thus $\ell_a$ is open and closed in $F^{-1}(\R)$, and being connected it is the component of $F^{-1}(\R)$ containing $a$.

(b) By the local description just used, every point of $\D$ has a neighborhood that meets at most one component of $F^{-1}(\R)$. Hence every compact subset of $\D$ meets only finitely many components, and an exhaustion of $\D$ by compact sets shows that there are countably many.

(c) Let $\ell$ be a component of $F^{-1}(\R)$. Its image $F(\ell)$ is a connected subset of $\R\setminus E$ and hence lies in a single component $I$ of $\R\setminus E$. The restriction of $F$ to $F^{-1}(I)$ is a covering map onto $I$, and $\ell$ is a component of $F^{-1}(I)$. Since the interval $I$ is simply connected, $F$ maps $\ell$ homeomorphically onto $I$. 
\end{proof}

\begin{lemma}\label{lem:separation}
Suppose that $M=\|\log F'\|_{\mathcal B}<\infty$, and let $\varepsilon_0$ and $\rho_0$ be as in Lemma~\ref{lem:injectivity}. Then any two distinct components $\ell$ and $\ell'$ of $F^{-1}(\R)$ satisfy
\begin{equation}\label{eq:separation}
d_\D(\ell,\ell')\ge\delta,\qquad \delta=\frac{\rho_0}{2}\ge\frac{\varepsilon_0}{4}=\frac{1}{16(1+M)} .
\end{equation}
\end{lemma}

\begin{proof}
Let $p\in\ell$. By Lemma~\ref{lem:reflection}(a), $F(R_{\ell'}p)=\overline{F(p)}=F(p)$. Since $p\notin\ell'$, the point $R_{\ell'}p$ is different from $p$, and the injectivity of $F$ on $D_\D(p,\rho_0)$ yields $d_\D(p,R_{\ell'}p)\ge\rho_0$. For $q\in\ell'$, the reflection $R_{\ell'}$ is an isometry fixing $q$, and so
\[
\rho_0\le d_\D(p,R_{\ell'}p)\le d_\D(p,q)+d_\D(q,R_{\ell'}p)=2d_\D(p,q).
\]
Taking the infimum over $p\in\ell$ and $q\in\ell'$ gives \eqref{eq:separation}.
\end{proof}

In particular, if $\log F'\in\mathcal B(\D)$, two distinct components of $F^{-1}(\R)$ cannot share an endpoint on $\T$.

\section{The length of the inverse image of the real line}\label{sec:length}

\subsection{Curvature and geodesics}

\begin{lemma}\label{lem:curvature}
Let $\ell$ be a component of $F^{-1}(\R)$, parametrized by Euclidean arclength $z=z(s)$, let $\kappa(s)$ be its signed curvature and set $u(s)=F(z(s))$. Then
\begin{equation}\label{eq:curvature-identity}
\frac{F''}{F'}\bigl(z(s)\bigr)\,z'(s)=\frac{d}{ds}\log|u'(s)|-i\kappa(s).
\end{equation}
In particular, if $\ell$ is a circular arc of Euclidean radius $r$, then $|F''/F'|\ge1/r$ on $\ell$.
\end{lemma}

\begin{proof}
The function $u$ is real valued and $u'(s)=F'(z(s))z'(s)\ne0$, so that $u''/u'=(d/ds)\log|u'|$ is real. By the chain rule,
\[
\frac{u''(s)}{u'(s)}=\frac{F''}{F'}\bigl(z(s)\bigr)z'(s)+\frac{z''(s)}{z'(s)},
\]
and $z''/z'=i\kappa$ for a curve parametrized by arclength. This gives \eqref{eq:curvature-identity}. Since $|z'|=1$, its imaginary part yields $|F''/F'|\ge|\kappa|=1/r$ on a circular arc of radius $r$.
\end{proof}

%We record the following elementary geometric formulas and estimates, which follow by direct calculation.
We need the following elementary geometric formulas and estimates, which follow by direct calculation.

\begin{lemma}\label{lem:geodesic}
Let $\ell$ be a complete hyperbolic geodesic of $\D$ that is a circular arc of Euclidean radius $r$. Let $a_\ell$ be its point nearest to $0$, and set $d_\ell=1-|a_\ell|$. Then
\[
|a_\ell|=\sqrt{1+r^2}-r,\qquad
\length(\ell)=2r\arctan\frac1r.
\]
Consequently,
\begin{alignat}{2}
\frac r2&\le d_\ell\le r,&\qquad \length(\ell)&\le\pi r\qquad\text{if }0<r\le1,\label{eq:small-geodesic}\\
d_\ell&>\frac12,&\qquad \length(\ell)&\le2\qquad\text{if }r>1.\label{eq:large-geodesic}
\end{alignat}
If $\ell$ is a diameter, we put $a_\ell=0$ and $d_\ell=1$; then $\length(\ell)=2$.
\end{lemma}

\subsection{The length estimate}

\begin{proof}[Proof of Theorem~\ref{thm:length}]
Write $M=\|\log F'\|_{\mathcal B}$ and
\[
Q=\int_{\D}\Bigl|\frac{F''(z)}{F'(z)}\Bigr|^2(1-|z|^2)\,dxdy.
\]
We may assume that $M$ and $Q$ are finite. Let $\varepsilon_0$ and $\delta$ be as in Lemmas~\ref{lem:injectivity} and~\ref{lem:separation}, and set
\[
\eta=\frac{\varepsilon_0}{64}=\frac{1}{256(1+M)} .
\]
For each component $\ell$ of $F^{-1}(\R)$ let $a_\ell$ and $d_\ell$ be as in Lemma~\ref{lem:geodesic}, let $r_\ell\in(0,\infty]$ be its Euclidean radius, with $r_\ell=\infty$ for a diameter, and let $D_\ell=D(a_\ell,\eta d_\ell)$. By Lemma~\ref{lem:reflection}(b), $F^{-1}(\R)$ has at most countably many components, and its length is the sum of their lengths.

We first show that the disks $D_\ell$ are pairwise disjoint. If $z\in D_\ell$, every point $w$ of the segment $[a_\ell,z]$ satisfies $1-|w|^2\ge(1-\eta)d_\ell$, while the segment has length less than $\eta d_\ell$. Integrating the hyperbolic metric along the segment, we see that
\[
d_\D(a_\ell,z)\le\frac{2\eta}{1-\eta}<\frac{\delta}{3},
\]
where the last inequality follows from $\delta\ge\varepsilon_0/4=16\eta$. Thus $D_\ell\subset D_\D(a_\ell,\delta/3)$. If $\ell\ne\ell'$, then $d_\D(a_\ell,a_{\ell'})\ge\delta$ by Lemma~\ref{lem:separation}, so $D_\ell\cap D_{\ell'}=\emptyset$. Note also that $D_\ell\subset\D$.

Next consider the components with $r_\ell\le1$. Let $\mu=|F''/F'|^2(1-|z|^2)\,dxdy$. Every $z\in D_\ell$ satisfies $1-|z|^2\ge(1-\eta)d_\ell\ge\frac34d_\ell$. Since $|F''/F'|^2$ is subharmonic, the mean value inequality gives
\[
\mu(D_\ell)\ge\frac34d_\ell\int_{D_\ell}\Bigl|\frac{F''(z)}{F'(z)}\Bigr|^2dxdy\ge\frac{3\pi}{4}\eta^2d_\ell^3\Bigl|\frac{F''}{F'}(a_\ell)\Bigr|^2 .
\]
By Lemma~\ref{lem:curvature} and \eqref{eq:small-geodesic}, the right-hand side is at least $\frac{3\pi}{4}\eta^2d_\ell^3r_\ell^{-2}\ge\frac{3\pi}{32}\eta^2r_\ell$, and therefore
\[
\length(\ell)\le\pi r_\ell\le\frac{32}{3}\eta^{-2}\mu(D_\ell).
\]
Summing over these components and using the disjointness of the disks, we find
\begin{equation}\label{eq:sum-small}
\sum_{r_\ell\le1}\length(\ell)\le\frac{32}{3}\eta^{-2}Q .
\end{equation}

It remains to consider the components with $r_\ell>1$, including a possible diameter. Each of them has $d_\ell>1/2$, by \eqref{eq:large-geodesic} and the convention $d_\ell=1$ for a diameter, so that the area of $D_\ell$ exceeds $\pi\eta^2/4$. These disks are disjoint and contained in $\D$, hence there are at most $4\eta^{-2}$ such components. Each of them has length at most $2$, and their total length is at most $8\eta^{-2}$. Together with \eqref{eq:sum-small} this gives
\[
\length(F^{-1}(\R))\le\eta^{-2}\Bigl(8+\frac{32}{3}Q\Bigr)\le C(1+M)^2(1+Q).\qedhere
\]
\end{proof}

\subsection{Uniform bounds and Carleson measures}

\begin{proposition}\label{prop:uniform-real}
If $\log F'\in\BMOA(\D)$, then
\[
\sup_{T\in\Aut(\D)}\length\bigl((F\circ T)^{-1}(\R)\bigr)\le C\bigl(1+\|\log F'\|_*\bigr)^4
\]
with an absolute constant $C$.
\end{proposition}

\begin{proof}
Let $T\in\Aut(\D)$ and set $B_0=\|\log F'\|_*+4\pi/\sqrt6$. Since, up to an additive constant, $\log(F\circ T)'=\log F'\circ T+\log T'$, Lemmas~\ref{lem:bmo-basic}(a) and~\ref{lem:log-aut} give $\|\log(F\circ T)'\|_*\le B_0$. Hence, by Lemma~\ref{lem:bmo-basic}(b),(c),
\[
\|\log(F\circ T)'\|_{\mathcal B}\le B_0,
\qquad
\int_{\D}\Bigl|\frac{(F\circ T)''}{(F\circ T)'}\Bigr|^2
(1-|z|^2)\,dxdy\le\pi B_0^2 .
\]
Theorem~\ref{thm:length}, applied to $F\circ T$, now gives
\[
\length\bigl((F\circ T)^{-1}(\R)\bigr)
\le C_0(1+B_0)^2(1+\pi B_0^2)
\le C\bigl(1+\|\log F'\|_*\bigr)^4,
\]
and the right-hand side does not depend on $T$.
\end{proof}

\begin{corollary}\label{cor:carleson-real}
For every universal covering $F$ of $\Omega$,
\begin{equation}\label{eq:S-length}
S\bigl(ds|_{F^{-1}(\R)}\bigr)=\sup_{T\in\Aut(\D)}\length\bigl((F\circ T)^{-1}(\R)\bigr),
\end{equation}
where $ds|_{F^{-1}(\R)}$ is arclength on $F^{-1}(\R)$. Consequently, arclength on $F^{-1}(\R)$ is a Carleson measure on $\D$ if and only if the right-hand side of \eqref{eq:S-length} is finite. This condition does not depend on the choice of $F$. If $\log F'\in\BMOA(\D)$, then $\|ds|_{F^{-1}(\R)}\|_{\CM}\le C(1+\|\log F'\|_*)^4$.
\end{corollary}

\begin{proof}
For $T\in\Aut(\D)$ we have $(F\circ T)^{-1}(\R)=T^{-1}(F^{-1}(\R))$, and
\[
\length\bigl(T^{-1}(F^{-1}(\R))\bigr)=\int_{F^{-1}(\R)}|(T^{-1})'(z)|\,ds(z)=\int_{F^{-1}(\R)}\frac{1-|a|^2}{|1-\overline az|^2}\,ds(z),\qquad a=T(0).
\]
Taking the supremum over $T$ gives \eqref{eq:S-length}. The remaining assertions follow from Lemma~\ref{lem:kernel}, from $(F\circ T)^{-1}(\R)=T^{-1}(F^{-1}(\R))$, and from Proposition~\ref{prop:uniform-real}.
\end{proof}

\section{Circles, lines and the Hayman--Wu property}\label{sec:circles}

\subsection{The components of \texorpdfstring{$\D\setminus F^{-1}(\R)$}{the complement of the inverse image of the real line}}

\begin{lemma}\label{lem:components}
Let $\{U_j\}$ be the components of $\D\setminus F^{-1}(\R)$.
\begin{enumerate}[label=\textup{(\alph*)}]
\item $F$ maps each $U_j$ conformally onto the upper half-plane $\Hp$ or onto the lower half-plane $\Hm$.
\item Each $U_j$ is hyperbolically convex, and $\partial U_j$ is a Jordan curve.
\item Every $\partial U_j$ is rectifiable and satisfies $\length(\partial U_j)\le\pi^2$. If $\length(F^{-1}(\R))<\infty$, then
\begin{equation}\label{eq:perimeter-sum}
\sum_j\length(\partial U_j)\le2\length(F^{-1}(\R))+2\pi .
\end{equation}
\end{enumerate}
\end{lemma}

\begin{proof}
(a) The set $\D\setminus F^{-1}(\R)$ is the disjoint union of $F^{-1}(\Hp)$ and $F^{-1}(\Hm)$. The restriction of $F$ to a component of $F^{-1}(\Hp)$ is a covering map onto the simply connected domain $\Hp$, hence a conformal map, and the same holds for $\Hm$.

(b) A component $\ell$ of $F^{-1}(\R)$ divides $\D$ into two open hyperbolic half-planes, and the connected set $U_j$ lies in one of them. If $p,q\in U_j$, the geodesic segment $[p,q]$ lies in the same half-plane for every $\ell$. Hence it does not meet $F^{-1}(\R)$, and being connected it lies in $U_j$. To see that $\partial U_j$ is a Jordan curve, consider the map $\mathcal K(z)=2z/(1+|z|^2)$ from the Poincar\'e model to the Klein model. It extends to a homeomorphism of $\overline\D$ and maps geodesic segments onto Euclidean segments. Thus $\mathcal K(U_j)$ is a bounded convex domain, whose boundary is a Jordan curve, and so is $\partial U_j=\mathcal K^{-1}(\partial\mathcal K(U_j))$.

(c) Every boundary point of $U_j$ in $\D$ belongs to $F^{-1}(\R)$. 
%Let $\ell$ be a component of $F^{-1}(\R)$. 
%By local finiteness, each $p\in\ell$ has a neighborhood $V$ that meets no other component of $F^{-1}(\R)$, so that the two sides of $V\setminus\ell$ lie in two components of $\D\setminus F^{-1}(\R)$. These components are distinct because $\ell$ separates $\D$. The local description shows that they remain unchanged as $p$ varies in a small subarc of $\ell$. Since $\ell$ is connected, the same two components are adjacent to it everywhere. Hence $\ell$ lies on the boundary of exactly two of the domains $U_j$, and a component of $F^{-1}(\R)$ that meets $\partial U_j$ is contained in $\partial U_j$.
For each component $\ell$ of $F^{-1}(\R)$ contained in $\partial U_j$, let $B_\ell$ be the open arc of $\T$ cut off by $\ell$ on the side opposite to $U_j$. The domain $U_j$ is the intersection of the open hyperbolic half-planes bounded by these components and containing $U_j$. 
%Indeed, if $p\in\D$ belongs to all these half-planes and $q\in U_j$, the geodesic segment from $q$ to $p$ lies in each of them. If this segment left $U_j$, its first exit point would belong to a component of $F^{-1}(\R)$ contained in $\partial U_j$, which is impossible, since the segment lies in the open half-plane bounded by that component. 
Since $U_j$ lies between any two of these components, the arcs $B_\ell$ are pairwise disjoint. Moreover, $\partial U_j\cap\T=\T\setminus\bigcup_\ell B_\ell$, because for $\zeta\in\T\setminus\bigcup_\ell B_\ell$ the geodesic ray from a point of $U_j$ to $\zeta$ stays in each of the half-planes. Let $\alpha_\ell,\beta_\ell\in\T$ be the endpoints of $\ell$. Since $\ell$ is a segment or a circular arc of angle at most $\pi$, we have
\[
\length(\ell)\le\frac\pi2|\alpha_\ell-\beta_\ell|\le\frac\pi2|B_\ell|.
\]
The pairwise disjoint arcs $B_\ell$ and their complement $\partial U_j\cap\T$ partition $\T$. Therefore
\[
\length(\partial U_j)
=|\partial U_j\cap\T|+\sum_\ell\length(\ell)
\le\frac\pi2\Bigl(|\partial U_j\cap\T|+\sum_\ell|B_\ell|\Bigr)
=\pi^2.
\]
Since a Jordan curve of finite one-dimensional Hausdorff measure is rectifiable \cite[Lemma~3.2]{Fa}, the curve $\partial U_j$ is rectifiable.

Every component of $F^{-1}(\R)$ lies on the boundary of exactly two of the $U_j$, and therefore
\[
\sum_j\length(\partial U_j\cap\D)=2\length(F^{-1}(\R)).
\]
If $U_i\ne U_j$, some component $\ell$ of $F^{-1}(\R)$ separates them. 
%Otherwise points $p\in U_i$ and $q\in U_j$ would lie on the same side of every component, and as in (b) the segment $[p,q]$ would not meet $F^{-1}(\R)$. 
Consequently $\overline{U_i}\cap\overline{U_j}\cap\T$ is contained in the set of the two endpoints of $\ell$. The family $\{U_j\}$ is countable, so the sets $\overline{U_j}\cap\T$ have pairwise intersections of measure zero, and the sum of their lengths is at most $2\pi$. This proves \eqref{eq:perimeter-sum}.
\end{proof}

\subsection{Circles and lines}

Let $\chi$ be the chordal metric on $\Cbar$,
\[
\chi(z,w)=\frac{2|z-w|}{\sqrt{1+|z|^2}\sqrt{1+|w|^2}},\qquad \chi(z,\infty)=\frac{2}{\sqrt{1+|z|^2}} ,
\]
let $B_\chi(p,r)$ be the open ball of center $p$ and radius $r$ for $\chi$, and let $\mathcal H^1_\chi$ be the one-dimensional Hausdorff measure with respect to $\chi$. Stereographic projection identifies $(\Cbar,\chi)$ with the unit sphere $S^2\subset\R^3$ with the Euclidean metric of $\R^3$.

\begin{lemma}\label{lem:chordal}
Let $\mathcal C$ be a circle or a straight line in $\C$, with the point $\infty$ added in the case of a line.
\begin{enumerate}[label=\textup{(\alph*)}]
\item $\mathcal H^1_\chi(\mathcal C\cap B_\chi(p,r))\le2\pi r$ for all $p\in\Cbar$ and $r>0$.
\item The Cayley maps $C_+(w)=(w-i)/(w+i)$ and $C_-(w)=(w+i)/(w-i)$ are isometries of $\chi$, and they map $\Hp$ and $\Hm$, respectively, onto $\D$.
\item Arclength on every Borel subset of $C_\pm(\mathcal C)\cap\D$ is a Carleson measure on $\D$ with Carleson norm at most $4\pi$.
\end{enumerate}
\end{lemma}

\begin{proof}
(a) Under inverse stereographic projection, $\mathcal C$ corresponds to a Euclidean circle $\mathcal S\subset S^2$ of radius $\varrho$, and $B_\chi(p,r)$ corresponds to a Euclidean ball of radius $r$ intersected with $S^2$. If $\varrho\le r$, the total length of $\mathcal S$ is at most $2\pi r$. Otherwise, if the intersection is nonempty, choose a point $a$ in it. Every other point in it is within Euclidean distance $2r$ of $a$, so the intersection lies in an arc of length at most $4\varrho\arcsin(r/\varrho)\le2\pi r$.

(b) A direct computation gives $|C_+(w)-C_+(w')|=2|w-w'|/(|w+i||w'+i|)$ and $1+|C_+(w)|^2=2(1+|w|^2)/|w+i|^2$, and the identity $\chi(C_+(w),C_+(w'))=\chi(w,w')$ follows. Since $C_-(w)=\overline{C_+(\overline w)}$ and complex conjugation is an isometry of $\chi$, the same holds for $C_-$.

(c) Let $X$ be a Borel subset of $C_+(\mathcal C)\cap\D$. For $u,v\in\overline\D$ we have $|u-v|\le\chi(u,v)\le2|u-v|$. Hence Euclidean length is at most $\mathcal H^1_\chi$ on subsets of $\overline\D$, and $D(\zeta,r)\cap\overline\D\subset B_\chi(\zeta,2r)$ for $\zeta\in\T$. For $\zeta\in\T$ and $r>0$, parts (a) and (b) yield
\[
\length\bigl(X\cap D(\zeta,r)\bigr)\le\mathcal H^1_\chi\bigl(C_+(\mathcal C)\cap B_\chi(\zeta,2r)\bigr)=\mathcal H^1_\chi\bigl(\mathcal C\cap B_\chi(C_+^{-1}(\zeta),2r)\bigr)\le4\pi r .
\]
The case of $C_-$ is the same.
\end{proof}

\begin{proposition}\label{prop:regular-sets}
Let $\mathcal C$ be a circle or a line in $\C$, and let $C_1$ be the constant in \eqref{eq:carleson-embedding} for $p=1$. If $\length(F^{-1}(\R))<\infty$, then
\[
\length\bigl(F^{-1}(\mathcal C)\bigr)\le(1+4C_1)\length\bigl(F^{-1}(\R)\bigr)+4\pi C_1 .
\]
\end{proposition}

\begin{proof}
Let $U_j$ be a component of $\D\setminus F^{-1}(\R)$ with $F(U_j)=\Hp$, let $g_j\colon\Hp\to U_j$ be the inverse of $F|_{U_j}$, and set $q_j=g_j\circ C_+^{-1}\colon\D\to U_j$. By Lemma~\ref{lem:components}, $\partial U_j$ is a rectifiable Jordan curve, and \eqref{eq:rectifiable} gives $q_j'\in H^1(\D)$ with $2\pi\|q_j'\|_{H^1}=\length(\partial U_j)$. Let $P_j=C_+(\mathcal C\cap\Hp)$. By Lemma~\ref{lem:chordal}(c) and \eqref{eq:carleson-embedding},
\[
\length\bigl(g_j(\mathcal C\cap\Hp)\bigr)=\int_{P_j}|q_j'|\,ds\le C_1\cdot4\pi\cdot\|q_j'\|_{H^1}=2C_1\length(\partial U_j).
\]
%The equality holds because $q_j$ is injective and conformal, so that $\length(q_j(X))=\int_X|q_j'|\,ds$ for every Borel set $X\subset\D$. 
The components mapped onto $\Hm$ are treated in the same way with $C_-$. Since $F^{-1}(\mathcal C)\cap U_j$ is $g_j(\mathcal C\cap\Hp)$ or the analogous set for $\Hm$, summation over $j$ and \eqref{eq:perimeter-sum} show that the part of $F^{-1}(\mathcal C)$ outside $F^{-1}(\R)$ has length at most $2C_1(2\length(F^{-1}(\R))+2\pi)$. The remaining part $F^{-1}(\mathcal C\cap\R)$ is contained in $F^{-1}(\R)$.
\end{proof}

\begin{proof}[Proof of Theorem~\ref{thm:uniform}]
By Proposition~\ref{prop:uniform-real}, the length of $(F\circ T)^{-1}(\R)$ is at most $C(1+\|\log F'\|_*)^4$, and Proposition~\ref{prop:regular-sets}, applied to the universal covering $F\circ T$, gives \eqref{eq:uniform-circles}. The statement about the Carleson norm is part of Corollary~\ref{cor:carleson-real}.
\end{proof}

\subsection{The Hayman--Wu property}\label{subsec:HW}

The Hayman--Wu property was defined in \eqref{eq:HW-intro}, where the constant is uniform over all universal coverings of $V$ and all lines $L$. Fern\'andez's theorem concerns complements of compact subsets of $\R$. The following lemma shows that $E$ is Carleson-homogeneous if its image under a suitable M\"obius transformation is a Carleson-homogeneous compact subset of $\R$.

\begin{lemma}\label{lem:transfer}
Let $E$ satisfy \eqref{eq:exterior}, let $a\in\R\setminus E$, and set $\Phi(w)=1/(w-a)$ and $K=\Phi(E)$. Then $K$ is a compact subset of $\R$ with $\diam K\ge2$. If $K$ is Carleson-homogeneous in the sense of \eqref{eq:compact-homogeneity}, then $E$ is Carleson-homogeneous.
\end{lemma}

\begin{proof}
Since $a\notin E$, the set $K$ is a compact subset of $\R$. It contains $0=\Phi(\infty)$, $\Phi(1)$ and $\Phi(-1)$, and $a\in(-1,1)$ gives
\[
\diam K\ge \Phi(1)-\Phi(-1)=\frac{2}{1-a^2}\ge2 .
\]
Let $c_K>0$ be a constant for which \eqref{eq:compact-homogeneity} holds for $K$, and set $d=\dist(a,E)>0$ and $r_0=\min\{d/2,1/2\}$. Let $x\in E\cap[-1,1]$, $0<r\le r_0$ and $I=(x-r,x+r)$. For $t\in I$ we have $d/2\le|t-a|<3$, and hence
\[
\frac19\le|\Phi'(t)|\le\frac{4}{d^2}\qquad(t\in I).
\]
The map $\Phi$ is monotone on $I$, and the lower bound for $|\Phi'|$ shows that $\Phi(I)$ contains $I'=(\Phi(x)-r/9,\Phi(x)+r/9)$. The scale $r/9$ is admissible in \eqref{eq:compact-homogeneity}, so $|K\cap I'|\ge c_Kr/9$. Since $I'\subset\Phi(I)$ and $\Phi$ is injective, $K\cap I'\subset \Phi(E\cap I)$. The upper bound for $|\Phi'|$ then gives
\[
|E\cap I|\ge\frac{d^2}{4}|K\cap I'|\ge\frac{c_Kd^2}{36}\,r .
\]
Put $c'=c_Kd^2/36$. If $r_0<r\le2$, monotonicity gives $|E\cap(x-r,x+r)|\ge c'r_0\ge\frac{c'r_0}{2}r$. If $r>2$, then $\R\setminus E\subset(-1,1)$ gives $|E\cap(x-r,x+r)|\ge2r-2\ge r$. Finally, if $x\ge1$, then $(x,x+r)\subset E$, and if $x\le-1$, then $(x-r,x)\subset E$. Thus \eqref{eq:homogeneity} holds with $C=\min\{1,c'r_0/2\}$.
\end{proof}

\begin{proof}[Proof of the equivalence of \textup{(i)--(iv)} in Theorem~\ref{thm:main}]
The implication (i)$\Rightarrow$(ii) is \cite[Theorem~1.1]{HZ}, and (ii)$\Rightarrow$(iii) is Corollary~\ref{cor:carleson-real}. By the same corollary, (iii) is equivalent to the finiteness of $\Lambda=\sup_{T\in\Aut(\D)}\length((F\circ T)^{-1}(\R))$.
If $\Lambda<\infty$, Proposition~\ref{prop:regular-sets}, applied to $F\circ T$, gives
\begin{equation}\label{eq:all-circles}
\length\bigl((F\circ T)^{-1}(\mathcal C)\bigr)
\le(1+4C_1)\Lambda+4\pi C_1
\end{equation}
for every $T\in\Aut(\D)$ and every circle or line $\mathcal C$. Since every universal covering of $\Omega$ has the form $F\circ T$, this proves (iii)$\Rightarrow$(iv). Conversely, (iv) implies $\Lambda<\infty$ by taking $L=\R$.

Assume (iv), choose $a\in\R\setminus E$, and let $\Phi$ and $K$ be as in Lemma~\ref{lem:transfer}. Then $V=\Phi(\Omega)=\Cbar\setminus K$, and $\widetilde F=\Phi\circ F$ is a meromorphic universal covering of $V$. For a line $L\subset\C$, the set $\mathcal C_L=\Phi^{-1}(L\cup\{\infty\})\cap\C$ is a circle or a line, and
\[
(\widetilde F\circ T)^{-1}(L\cap V)
\subset (F\circ T)^{-1}(\mathcal C_L).
\]
Every universal covering of $V$ is of the form $\widetilde F\circ T$, and \eqref{eq:all-circles} bounds these lengths uniformly in $T$ and $L$. Thus $V$ has the Hayman--Wu property. Fern\'andez's theorem implies that $K$ is Carleson-homogeneous, and Lemma~\ref{lem:transfer} gives (i).
\end{proof}

\begin{remark}
The M\"obius transformation $\Phi$ is used only after the length estimates have been established for $F$, since $\Phi\circ F$ has poles and the estimates of Section~\ref{sec:length} do not apply to it.
\end{remark}

\section{The Schwarzian derivative}\label{sec:schwarzian}

In this section $f$ is locally univalent and analytic in $\D$, $S_f=(f''/f')'-\frac12(f''/f')^2$ is
the Schwarzian derivative of $f$, and
\[
\nu_f=(1-|z|^2)^3|S_f(z)|^2\,dxdy .
\]
We give a direct proof that $\log f'\in\BMOA(\D)$ if and only if $\nu_f$ is a Carleson measure. The equivalence also follows from Gr\"ohn, Nicolau and R\"atty\"a \cite[Theorem~4(i)]{GNR}, applied to $u=(f')^{-1/2}$, which satisfies $u''+\frac12S_fu=0$. For conformal maps it goes back to \cite{AZ}, and universal coverings of uniformly perfect domains are discussed in \cite[Section~4]{HZ}. The proof below first assumes a Bloch bound and then derives it from the Schwarzian condition.

\subsection{Transformation rules}

We shall use the following standard transformation rules for the pre-Schwarzian and Schwarzian derivatives.

\begin{lemma}\label{lem:transformation}
Let $T\in\Aut(\D)$, $z\in\D$ and $w=T(z)$. Then
\begin{equation}\label{eq:invariant-preschwarzian}
(1-|z|^2)\frac{(f\circ T)''}{(f\circ T)'}(z)-2\overline z=\Bigl((1-|w|^2)\frac{f''}{f'}(w)-2\overline w\Bigr)\frac{T'(z)}{|T'(z)|}.
\end{equation}
In particular $\|\log(f\circ T)'\|_{\mathcal B}\le\|\log f'\|_{\mathcal B}+4$. Moreover $S_{f\circ T}=(S_f\circ T)(T')^2$, so that 
\[
(1-|z|^2)^2|S_{f\circ T}(z)|=(1-|w|^2)^2|S_f(w)|.
\]
\end{lemma}

\begin{proof}
The identity \eqref{eq:invariant-preschwarzian} follows by substituting
\[
\frac{(f\circ T)''}{(f\circ T)'}
=\Bigl(\frac{f''}{f'}\circ T\Bigr)T'+\frac{T''}{T'}
\]
and using the identities $(1-|z|^2)|T'(z)|=1-|T(z)|^2$
and
\[
(1-|z|^2)\frac{T''(z)}{T'(z)}-2\overline z
=-2\overline{T(z)}\,\frac{T'(z)}{|T'(z)|},
\]
which hold for every automorphism of $\D$. The estimate for the Bloch seminorm follows from \eqref{eq:invariant-preschwarzian}, since $|z|<1$ and $|w|<1$. The formula for the Schwarzian derivative is the chain rule, since $S_T=0$.
\end{proof}

For the universal covering $F$, the chain rule shows that $S_F=(S_F\circ\gamma)(\gamma')^2$ for $\gamma\in G$, so that $(1-|z|^2)^2|S_F(z)|$ is invariant under $G$.

The substitution $z=\sigma_a(v)$ will be used repeatedly. Its Jacobian is $|\sigma_a'(v)|^2$, and we have $1-|z|^2=(1-|v|^2)|\sigma_a'(v)|$ and $|\sigma_a'(z)|=1/|\sigma_a'(v)|$. Hence, for every nonnegative Borel function $h$ on $\D$,
\begin{equation}\label{eq:change-of-variables}
\int_{\D}h(z)\frac{1-|a|^2}{|1-\overline az|^2}\,dxdy=\int_{\D}h(\sigma_a(z))\,|\sigma_a'(z)|\,dxdy .
\end{equation}
We also need the analogue of \eqref{eq:area-identity} for second derivatives. If $\varphi(z)=\sum_{n\ge0}c_nz^n$, then
\begin{equation}\label{eq:second-derivative}
\int_{\D}|\varphi''(z)|^2(1-|z|^2)^3\,dxdy=6\pi\sum_{n\ge2}\frac{n(n-1)}{(n+1)(n+2)}|c_n|^2\le6\pi\|\varphi-\varphi(0)\|_{H^2}^2 ,
\end{equation}
which follows from $\int_0^1s^{2n-3}(1-s^2)^3\,ds=3/\bigl((n-1)n(n+1)(n+2)\bigr)$ for $n\ge2$.

\subsection{From BMOA to the Schwarzian derivative}

\begin{proposition}\label{prop:bmoa-schwarzian}
Suppose that $\varphi=\log f'\in\BMOA(\D)$, and let $B=\|\varphi\|_*$. 
For $d\mu_1=(1-|z|^2)|\varphi'|^2\,dxdy$, $d\mu_2=(1-|z|^2)^3|\varphi''|^2\,dxdy$, and $d\nu_f$,
there is an absolute constant $C>0$ such that
\[
\|\mu_1\| \le CB^2,\quad \|\mu_2\| \le CB^2, \quad {\rm and} \quad \|\nu_f\| \le C(B^2+B^4).
\]
%\begin{align*}
%\bigl\|(1-|z|^2)|\varphi'(z)|^2\,dxdy\bigr\|_{\CM}&\le CB^2,\\
%\bigl\|(1-|z|^2)^3|\varphi''(z)|^2\,dxdy\bigr\|_{\CM}&\le CB^2,\\
%\|\nu_f\|_{\CM}&\le C(B^2+B^4).
%\end{align*}
\end{proposition}

\begin{proof}
%Write $d\mu_1=(1-|z|^2)|\varphi'|^2\,dxdy$ and $d\mu_2=(1-|z|^2)^3|\varphi''|^2\,dxdy$. 
By Lemma~\ref{lem:kernel}, it suffices to estimate $S(\mu_1)$, $S(\mu_2)$ and $S(\nu_f)$. Fix $a\in\D$ and set $\psi_a=\varphi\circ\sigma_a-\varphi(a)$, so that $\|\psi_a\|_{H^2}\le B$. The change of variables \eqref{eq:change-of-variables} and Lemma~\ref{lem:bmo-basic}(c) give
\[
\int_{\D}\frac{1-|a|^2}{|1-\overline az|^2}\,d\mu_1(z)
=\int_{\D}|\psi_a'(z)|^2(1-|z|^2)\,dxdy
\le\pi B^2.
\]
Thus $S(\mu_1)\le\pi B^2$.

For $\mu_2$, the same change of variables gives the integral of
$|\varphi''(\sigma_a(z))|^2|\sigma_a'(z)|^4(1-|z|^2)^3$.
Since
\[
(\varphi''\circ\sigma_a)(\sigma_a')^2
=\psi_a''-\psi_a'\frac{\sigma_a''}{\sigma_a'},
\qquad
(1-|z|^2)\left|\frac{\sigma_a''(z)}{\sigma_a'(z)}\right|\le4,
\]
this integral is at most
\[
C\int_{\D}\bigl(|\psi_a''(z)|^2(1-|z|^2)^3
+|\psi_a'(z)|^2(1-|z|^2)\bigr)\,dxdy
\le CB^2
\]
by \eqref{eq:second-derivative} and Lemma~\ref{lem:bmo-basic}(c). Taking the supremum over $a$ yields $S(\mu_2)\le CB^2$.

Finally, $|S_f|^2\le2|\varphi''|^2+\frac12|\varphi'|^4$, and the Bloch bound $(1-|z|^2)|\varphi'(z)|\le B$ gives $\nu_f\le2\mu_2+\frac12B^2\mu_1$. Hence
\[
S(\nu_f)\le2S(\mu_2)+\tfrac12B^2S(\mu_1)
\le C(B^2+B^4).
\]
Lemma~\ref{lem:kernel} proves all three estimates.
\end{proof}

\subsection{From the Schwarzian derivative to BMOA}

\begin{lemma}\label{lem:normal-family}
Let $M\ge0$, and let $g$ be a locally univalent analytic function in $\D$ such that
\[
\|\log g'\|_{\mathcal B}\le M
\qquad\text{and}\qquad
\left|\frac{g''(0)}{g'(0)}\right|\ge\frac52 .
\]
Then
\[
\sup_{|u|<1/2}(1-|u|^2)^2|S_g(u)|\ge\eta_M ,
\]
where $\eta_M>0$ is a constant depending only on $M$.
\end{lemma}

\begin{proof}
Suppose the contrary. Then there are functions $g_n$ with the stated properties such that $(1-|u|^2)^2|S_{g_n}(u)|\le1/n$ for $|u|<1/2$, and affine postcomposition allows us to assume that $g_n(0)=0$ and $g_n'(0)=1$. Choose the analytic branch with $\log g_n'(0)=0$. Integration of the Bloch bound along radii gives $|\log g_n'(u)|\le M\artanh|u|$. Hence $\{g_n'\}$ is locally bounded and locally bounded away from $0$, and a subsequence of $\{g_n\}$ converges locally uniformly to an analytic function $g$ with $g'\ne0$, $|g''(0)/g'(0)|\ge5/2$ and $S_g=0$ on $D(0,1/2)$. By analyticity $S_g=0$ in $\D$, so $g(u)=(au+b)/(cu+d)$ is a M\"obius transformation without poles in $\D$. Thus $d\ne0$ and $|c|\le|d|$, and $|g''(0)/g'(0)|=2|c/d|\le2$, a contradiction.
\end{proof}

\begin{lemma}\label{lem:main-estimate}
Let $g$ be a locally univalent analytic function in $\D$, and let $M,K\ge0$. Suppose that
\[
\|\log g'\|_{\mathcal B}\le M
\qquad\text{and}\qquad
\int_{\D}(1-|z|^2)^3|S_g(z)|^2\,dxdy\le K.
\]
Then
\[
\int_{\D}\left|\frac{g''(z)}{g'(z)}\right|^2(1-|z|^2)\,dxdy
\le C,
\]
where the constant $C=C(M,K)$ depends only on $M$ and $K$.
\end{lemma}

\begin{proof}
Assume first that $g''/g'$ is bounded in $\D$, so that the integral $I_1$ on the left-hand side is finite. Write $\varphi(z)=\log g'(z)-\log g'(0)=\sum_{n\ge1}c_nz^n$ and
\[
I_2=\int_{\D}|\varphi''(z)|^2(1-|z|^2)^3\,dxdy.
\]
Since $c_1=\varphi'(0)$, we have $|c_1|\le M$. For $n\ge2$, Cauchy's estimate for $\varphi'$ on $|z|=r_n$, where $r_n=1-1/n$, gives
\[
n|c_n|\le \frac{M}{(1-r_n^2)r_n^{\,n-1}}.
\]
Since $n(1-r_n^2)\ge1$ and $r_n^{\,n-1}\ge e^{-1}$, it follows that $|c_n|\le eM$. Let $N\ge2$ and set $\beta_N=(N+2)/(6(N-1))$. For $n\ge N$ we have
\[
\frac{n}{n+1}\le\beta_N\cdot\frac{6n(n-1)}{(n+1)(n+2)},
\]
and therefore \eqref{eq:area-identity} and \eqref{eq:second-derivative} give
\begin{equation}\label{eq:I1-I2}
I_1\le\beta_NI_2+\pi Ne^2M^2 .
\end{equation}

Next, $\varphi''=S_g+\frac12(\varphi')^2$, so that 
\[
|\varphi''|^2\le(1+\tau^{-1})|S_g|^2+\frac14(1+\tau)|\varphi'|^4
\] 
for every $\tau>0$. Let $X=\{z\in\D:(1-|z|^2)|\varphi'(z)|>9/2\}$. Outside $X$ we have 
\[
(1-|z|^2)^3|\varphi'(z)|^4\le\frac{81}{4}(1-|z|^2)|\varphi'(z)|^2, 
\]
and on $X$ we have $(1-|z|^2)^3|\varphi'(z)|^4\le M^4(1-|z|^2)^{-1}$. Consequently,
\begin{equation}\label{eq:I2-bound}
I_2\le(1+\tau^{-1})K+\frac{1+\tau}{4}\Bigl(\frac{81}{4}I_1+M^4\int_X\frac{dxdy}{1-|z|^2}\Bigr).
\end{equation}

We now estimate the integral over $X$. Let $\eta_*=\eta_{M+4}/2$, with $\eta_{M+4}$ given by Lemma~\ref{lem:normal-family}. For $v\in X$ let $h=g\circ\sigma_v$. By Lemma~\ref{lem:transformation}, $\|\log h'\|_{\mathcal B}\le M+4$. Since $\sigma_v'(0)=-(1-|v|^2)$ and $\sigma_v''(0)/\sigma_v'(0)=2\overline v$, the number $h''(0)/h'(0)=2\overline v-(1-|v|^2)\varphi'(v)$ has modulus larger than $5/2$. Lemma~\ref{lem:normal-family} gives a point $u$ with $|u|<1/2$ and $(1-|u|^2)^2|S_h(u)|\ge\eta_*$. The point $\xi_v=\sigma_v(u)$ then satisfies
\[
d_\D(v,\xi_v)<\log3,\qquad (1-|\xi_v|^2)^2|S_g(\xi_v)|\ge\eta_* .
\]
Put $R=1+\log3$, and choose a maximal, at most countable family $\{v_k\}\subset X$ with $d_\D(v_k,v_j)\ge2R$ for $k\ne j$. The disks $D_\D(v_k,R)$ are pairwise disjoint, and $X\subset\bigcup_kD_\D(v_k,2R)$ by maximality. The automorphism $\sigma_v$ maps $D(0,\tanh R)$ onto $D_\D(v,2R)$, and $|\sigma_v'(z)|\le(1-|v|^2)(1-\tanh R)^{-2}$ for $|z|<\tanh R$. Hence the substitution $z\mapsto\sigma_v(z)$ gives
\[
\int_{D_\D(v,2R)}\frac{dxdy}{1-|z|^2}=\int_{D(0,\tanh R)}\frac{|\sigma_v'(z)|}{1-|z|^2}\,dxdy\le C(1-|v|^2).
\]
On the other hand, let $t=\tanh(1/2)$ and $\xi=\xi_{v_k}$. The substitution $z\mapsto\sigma_\xi(z)$, the transformation rule for $S_{g\circ\sigma_\xi}$ and the subharmonicity of $|S_{g\circ\sigma_\xi}|^2$ give
\begin{align*}
\int_{D_\D(\xi,1)}(1-|z|^2)^3|S_g(z)|^2\,dxdy&=\int_{D(0,t)}(1-|z|^2)^3|S_{g\circ\sigma_\xi}(z)|^2|\sigma_\xi'(z)|\,dxdy\\
&\ge\frac{\pi t^2(1-t^2)^3}{4}\,(1-|\xi|^2)\eta_*^2 .
\end{align*}
Here we used $|\sigma_\xi'(z)|\ge(1-|\xi|^2)/4$ and $|S_{g\circ\sigma_\xi}(0)|=(1-|\xi|^2)^2|S_g(\xi)|$. Since $D_\D(\xi_{v_k},1)\subset D_\D(v_k,R)$, these disks are pairwise disjoint, and \eqref{eq:harnack} gives $1-|v_k|^2\le3(1-|\xi_{v_k}|^2)$. Therefore
\[
\int_X\frac{dxdy}{1-|z|^2}\le C\sum_k(1-|v_k|^2)\le\frac{C}{\eta_*^2}\sum_k\int_{D_\D(\xi_{v_k},1)}(1-|z|^2)^3|S_g(z)|^2\,dxdy\le\frac{C}{\eta_*^2}K .
\]

Choose $N=100$ and $\tau=1/10$. Then $\beta_N=17/99$ and $q=(1+\tau)\beta_N\cdot81/16=153/160<1$. Inserting \eqref{eq:I2-bound} and the last estimate into \eqref{eq:I1-I2}, we find $(1-q)I_1\le C(M,K)$, which proves the lemma when $g''/g'$ is bounded.

In general, let $0<r<1$ and $g_r(z)=g(rz)$. Then $g_r''/g_r'$ is bounded, and $\|\log g_r'\|_{\mathcal B}\le\|\log g'\|_{\mathcal B}$ because $(1-|z|^2)r\le1-r^2|z|^2$. Moreover $S_{g_r}(z)=r^2S_g(rz)$, and since $r^2-|z|^2\le r^2(1-|z|^2)$, the substitution $z\mapsto z/r$ gives
\[
\int_{\D}(1-|z|^2)^3|S_{g_r}(z)|^2\,dxdy=\int_{D(0,r)}\frac{(r^2-|z|^2)^3}{r^4}|S_g(z)|^2\,dxdy\le K.
\]
The first part of the proof applies to $g_r$, and Fatou's lemma completes the proof as $r\to1$.
\end{proof}

\begin{proposition}\label{prop:schwarzian-bmoa}
Suppose that $\|\log f'\|_{\mathcal B}\le M$ and that $\nu_f$ is a Carleson measure on $\D$. Then $\log f'\in\BMOA(\D)$, and $\|\log f'\|_*$ is bounded by a constant depending only on $M$ and $\|\nu_f\|_{\CM}$.
\end{proposition}

\begin{proof}
Fix $a\in\D$ and let $g=f\circ\sigma_a$. By Lemma~\ref{lem:transformation}, $\|\log g'\|_{\mathcal B}\le M+4$, and by \eqref{eq:change-of-variables} and Lemma~\ref{lem:kernel},
\[
\int_{\D}(1-|z|^2)^3|S_g(z)|^2\,dxdy=\int_{\D}\frac{1-|a|^2}{|1-\overline az|^2}\,d\nu_f(z)\le C_S\|\nu_f\|_{\CM}.
\]
Lemma~\ref{lem:main-estimate} bounds $\int_{\D}|g''(z)/g'(z)|^2(1-|z|^2)\,dxdy$ by a constant depending only on $M$ and $\|\nu_f\|_{\CM}$. Since $(\log f'\circ\sigma_a)'=g''/g'-\sigma_a''/\sigma_a'$, Lemma~\ref{lem:bmo-basic}(c) and the bound for $\|k_{\sigma_a}\|_{H^2}$ in Lemma~\ref{lem:log-aut} yield
\[
\frac{\pi}{2}\|\log f'\circ\sigma_a-\log f'(a)\|_{H^2}^2\le2\int_{\D}\Bigl|\frac{g''(z)}{g'(z)}\Bigr|^2(1-|z|^2)\,dxdy+2\pi\cdot\frac{2\pi^2}{3}.
\]
The right-hand side is bounded independently of $a$.
\end{proof}

The assumption $\|\log f'\|_{\mathcal B}\le M$ in Proposition~\ref{prop:schwarzian-bmoa} can be removed. We first note that the Carleson condition bounds the Schwarzian derivative pointwise.

\begin{lemma}\label{lem:schwarzian-bounded}
If $\nu_f$ is a Carleson measure, then
\begin{equation}\label{eq:schwarzian-bounded}
\sup_{z\in\D}(1-|z|^2)^2|S_f(z)|\le16\,\|\nu_f\|_{\CM}^{1/2}.
\end{equation}
\end{lemma}

\begin{proof}
Let $z\in\D$, set $\varrho=(1-|z|)/2$, and let $\zeta=z/|z|$ (any $\zeta\in\T$ if $z=0$). Then $D(z,\varrho)\subset D(\zeta,3\varrho)$, and $1-|u|^2\ge\varrho$ for $u\in D(z,\varrho)$. The subharmonicity of $|S_f|^2$ gives
\[
\pi\varrho^5|S_f(z)|^2\le\varrho^3\int_{D(z,\varrho)}|S_f|^2\,dxdy\le\nu_f\bigl(D(\zeta,3\varrho)\bigr)\le3\varrho\,\|\nu_f\|_{\CM}.
\]
Hence $(1-|z|^2)^2|S_f(z)|\le16\varrho^2|S_f(z)|\le16(3/\pi)^{1/2}\|\nu_f\|_{\CM}^{1/2}$.
\end{proof}

\begin{proposition}\label{prop:nehari}
Let $f$ be a locally univalent analytic function in $\D$ such that
\[
N=\sup_{z\in\D}(1-|z|^2)^2|S_f(z)|<\infty .
\]
Then $\|\log f'\|_{\mathcal B}\le8\sqrt{1+N}+2$.
\end{proposition}

\begin{proof}
Fix $a\in\D$, and set $g=f\circ\sigma_a$, $r=1/(2\sqrt{1+N})$ and $h(z)=g(rz)$. By Lemma~\ref{lem:transformation}, $(1-|z|^2)^2|S_g(z)|\le N$ in $\D$. Since $S_h(z)=r^2S_g(rz)$ and $1-|z|^2\le1-r^2|z|^2$, we have
\[
(1-|z|^2)^2|S_h(z)|\le r^2N\,\frac{(1-|z|^2)^2}{(1-r^2|z|^2)^2}\le r^2N<2\qquad(z\in\D).
\]
By Nehari's univalence criterion \cite{Ne}, $h$ is univalent in $\D$, and the Bieberbach inequality for the second coefficient \cite[Chapter~1]{Po75} gives $|h''(0)/h'(0)|\le4$. On the other hand, $\sigma_a'(0)=-(1-|a|^2)$ and $\sigma_a''(0)/\sigma_a'(0)=2\overline a$, so that
\[
\frac{h''(0)}{h'(0)}=r\,\frac{g''(0)}{g'(0)}=r\Bigl(2\overline a-(1-|a|^2)\frac{f''(a)}{f'(a)}\Bigr).
\]
Hence $(1-|a|^2)|f''(a)/f'(a)|\le4/r+2=8\sqrt{1+N}+2$, and the assertion follows by taking the supremum over $a\in\D$.
\end{proof}

\begin{theorem}\label{thm:schwarzian}
Let $f$ be a locally univalent analytic function in $\D$. Then $\log f'\in\BMOA(\D)$ if and only if $\nu_f\in\CM(\D)$.
\end{theorem}

\begin{proof}
If $\log f'\in\BMOA(\D)$, the Carleson condition holds by Proposition~\ref{prop:bmoa-schwarzian}. Conversely, if $\nu_f$ is a Carleson measure, Lemma~\ref{lem:schwarzian-bounded} and Proposition~\ref{prop:nehari} give $\log f'\in\mathcal B(\D)$, and Proposition~\ref{prop:schwarzian-bmoa} gives $\log f'\in\BMOA(\D)$.
\end{proof}

\begin{proof}[Completion of the proof of Theorem~\ref{thm:main}]
Conditions (i)--(iv) are equivalent by Section~\ref{subsec:HW}, and Theorem~\ref{thm:schwarzian} with $f=F$ gives (ii)$\Leftrightarrow$(v).
\end{proof}

In particular, condition (v) readily implies that $\Omega$ is uniformly perfect. Indeed, Lemma~\ref{lem:schwarzian-bounded} and Proposition~\ref{prop:nehari} give $\log F'\in\mathcal B(\D)$, and Pommerenke's theorem \cite{Po84} applies.

\subsection{The little Bloch space}

The little Bloch space $\mathcal B_0(\D)$ consists of analytic functions $\varphi$ with $(1-|z|^2)|\varphi'(z)|\to0$ as $|z|\to1$, and $\VMOA(\D)\subset\mathcal B_0(\D)$ \cite{Po78}. The corresponding vanishing conditions never hold for the coverings $F:\D \to \Omega=\Cbar \setminus E$ considered here.

\begin{proposition}\label{prop:little-bloch}
Neither $(1-|z|^2)|F''(z)/F'(z)|$ nor $(1-|z|^2)^2|S_F(z)|$ tends to $0$ as $|z|\to1$. In particular, $\log F'\notin\VMOA(\D)$.
\end{proposition}

\begin{proof}
Suppose first that $\Omega$ is not simply connected. Then the Fuchsian group $G$ is infinite. For $\gamma\in G$, the identity $F\circ\gamma=F$ gives $(F''/F')\circ\gamma\cdot\gamma'=F''/F'-\gamma''/\gamma'$, and hence
\[
(1-|\gamma(z)|^2)\Bigl|\frac{F''}{F'}(\gamma(z))\Bigr|=(1-|z|^2)\Bigl|\frac{F''}{F'}(z)-\frac{\gamma''}{\gamma'}(z)\Bigr| .
\]
Let $\gamma_n\in G$ be distinct. Then $|\gamma_n(z)|\to1$ for each $z\in\D$, since $G$ acts discontinuously. If $\log F'\in\mathcal B_0(\D)$, the left-hand side tends to $0$, so that $\gamma_n''/\gamma_n'\to F''/F'$ pointwise in $\D$. Writing $\gamma_n(z)=e^{i\theta_n}(z-\alpha_n)/(1-\overline\alpha_nz)$, we have $\gamma_n''/\gamma_n'(z)=2\overline\alpha_n/(1-\overline\alpha_nz)$ and $|\alpha_n|=|\gamma_n^{-1}(0)|\to1$. Passing to a subsequence with $\alpha_n\to\alpha\in\T$, we find $F''/F'=2\overline\alpha/(1-\overline\alpha z)$. Then $F'$ is a constant multiple of $(1-\overline\alpha z)^{-2}$, so $F$ is a M\"obius transformation and in particular injective, which contradicts $F\circ\gamma=F$ for $\gamma\ne\mathrm{id}$. The function $(1-|z|^2)^2|S_F(z)|$ is $G$-invariant, so if it tended to $0$ at $\T$, it would vanish identically. Then $S_F=0$ and $F$ would again be a M\"obius transformation.

Suppose now that $\Omega$ is simply connected. Its complement $E$ is connected in $\Cbar$, so $\Rbar\setminus E$ is a single interval $(\alpha,\beta)$. Thus $E=[-\infty,\alpha]\cup[\beta,+\infty]$ with $-1\le\alpha<\beta\le1$, and $F=A\circ F_0\circ T$ with $T\in\Aut(\D)$, the affine map $A(w)=\frac12(\beta-\alpha)w+\frac12(\alpha+\beta)$ and
\[
F_0(z)=\frac{2z}{1+z^2}.
\]
Indeed, $F_0$ is the reciprocal of the Joukowski map $z\mapsto\frac12(z+z^{-1})$, and it maps $\D$ conformally onto $\C\setminus\bigl((-\infty,-1]\cup[1,\infty)\bigr)$. A direct computation gives
\[
\frac{F_0''}{F_0'}(z)=-\frac{2z}{1-z^2}-\frac{4z}{1+z^2},\qquad S_{F_0}(z)=-\frac{6}{(1-z^2)^2},
\]
so that $(1-r^2)F_0''(r)/F_0'(r)\to-2$ as $r\to1$ and $(1-r^2)^2|S_{F_0}(r)|=6$ for $0<r<1$. The function $\log T'$ belongs to $\mathcal B_0(\D)$, and $\varphi\circ T\in\mathcal B_0(\D)$ if and only if $\varphi\in\mathcal B_0(\D)$. Since $\log F'=\log F_0'\circ T+\log T'$ up to an additive constant, it follows that $\log F'\notin\mathcal B_0(\D)$. Finally, $(1-|z|^2)^2|S_F(z)|=(1-|T(z)|^2)^2|S_{F_0}(T(z))|$ by Lemma~\ref{lem:transformation}, and this quantity equals $6$ at the points $z=T^{-1}(r)$, $0<r<1$.
\end{proof}

\section{Interpolating sequences}\label{sec:interpolation}

A sequence $W\subset\D$ is an interpolating sequence for $H^\infty$ if every bounded function on $W$ is the restriction of a function in $H^\infty$. By Carleson's interpolation theorem \cite[Chapter~VII]{Ga}, this holds if and only if $W$ is separated in the hyperbolic metric and $\sum_{w\in W}(1-|w|^2)\delta_w$ is a Carleson measure, and the interpolation constant depends only on the separation constant and the Carleson norm. Moreover, this is equivalent to the condition that $W$ is \emph{uniformly separated}, that is, $\inf_{w\in W}\prod_{v\in W\setminus\{w\}}\rho(w,v)>0$, where $\rho(u,v)=|\sigma_u(v)|$ is the pseudo-hyperbolic distance.

\subsection{Interpolation estimates}\label{subsec:interpolation-estimate}
The following density estimate will be used below.
This is a quantitative form of the fact that $\log F' \in {\mathcal B}(\D)$ implies that
$\Omega$ is uniformly perfect.

\begin{lemma}\label{lem:density}
If $M=\|\log F'\|_{\mathcal B}<\infty$, then
\[
\frac{1}{16(1+M)}\le\lambda_\Omega(w)\dist(w,E)\le2\qquad(w\in\Omega).
\]
\end{lemma}

\begin{proof}
The upper bound holds in every domain and follows by comparison with the disk centered at $w$ with radius $\dist(w,E)$, which is contained in $\Omega$. For the lower bound, let $w=F(a)$ and $d=1-|a|$, and let $\varepsilon_0$ be as in Lemma~\ref{lem:injectivity}. Since $F$ is injective on $D(a,\varepsilon_0d)$, the Koebe one-quarter theorem \cite[Chapter~1]{Po75} shows that $F(D(a,\varepsilon_0d))\subset\Omega$ contains the disk of center $w$ and radius $\varepsilon_0d|F'(a)|/4$. Therefore
\[
\dist(w,E)\ge\frac{\varepsilon_0d|F'(a)|}{4}\ge\frac{\varepsilon_0}{8}(1-|a|^2)|F'(a)|=\frac{\varepsilon_0}{4\lambda_\Omega(w)} ,
\]
which is the assertion.
\end{proof}

\begin{proposition}\label{prop:interpolation}
Suppose that $\log F'\in\BMOA(\D)$, and let $M=\|\log F'\|_{\mathcal B}$. Let $0<\kappa\le1$, and let $Z\subset\R\setminus E$ satisfy
\begin{equation}\label{eq:sparse}
|z-z'|\ge\kappa\max\{\dist(z,E),\dist(z',E)\}\qquad(z,z'\in Z,\ z\ne z').
\end{equation}
Set $W=F^{-1}(Z)$ and $\mu_W=\sum_{w\in W}(1-|w|^2)\delta_w$. Then
\begin{enumerate}[label=\textup{(\alph*)}]
\item $d_\D(w,w')\ge\kappa/(64(1+M))$ for all distinct $w,w'\in W$, and
\item $\|\mu_W\|_{\CM}\le512\,e^2(1+M)\kappa^{-1}\|ds|_{F^{-1}(\R)}\|_{\CM}\le C\kappa^{-1}(1+\|\log F'\|_*)^5$ with an absolute constant $C$.
\end{enumerate}
Consequently $W$ is an interpolating sequence for $H^\infty$, with interpolation constant depending only on $\kappa$ and $\|\log F'\|_*$.
\end{proposition}

\begin{proof}
For $z\in Z$ write $d_z=\dist(z,E)$ and $I_z=(z-\kappa d_z/3,z+\kappa d_z/3)$. These intervals are contained in $\R\setminus E$ and are pairwise disjoint, since a common point of $I_z$ and $I_{z'}$ would give $|z-z'|<\frac{2\kappa}{3}\max\{d_z,d_{z'}\}$. For $x\in I_z$ we have $\frac23d_z\le\dist(x,E)\le\frac43d_z$. By Lemma~\ref{lem:density},
\[
\frac{c_0}{\dist(x,E)}\le\lambda_\Omega(x)\le\frac{2}{\dist(x,E)},\qquad c_0=\frac{1}{16(1+M)} .
\]

(a) If $w$ and $w'$ lie on different components of $F^{-1}(\R)$, then $d_\D(w,w')\ge\delta\ge c_0$ by Lemma~\ref{lem:separation}. Suppose that they lie on the same component $\ell$. By Lemma~\ref{lem:reflection}(c), $z=F(w)$ and $z'=F(w')$ are distinct points of the same component of $\R\setminus E$, and the arc of $\ell$ between $w$ and $w'$ is mapped onto the segment $[z,z']$. Since $F$ is a local isometry and $\ell$ is a geodesic,
\[
d_\D(w,w')=\int_{[z,z']}\lambda_\Omega(x)\,dx\ge\int_{I_z\cap[z,z']}\frac{c_0\,dx}{\dist(x,E)}\ge\frac{3c_0}{4d_z}\cdot\frac{\kappa d_z}{3}=\frac{c_0\kappa}{4},
\]
because $|z-z'|\ge\kappa d_z$ implies that $[z,z']$ contains one half of $I_z$. This proves (a).

(b) Fix $w\in W$, let $z=F(w)$, let $\ell$ be the component of $F^{-1}(\R)$ through $w$, and let $\Gamma_w\subset\ell$ be the arc that $F$ maps onto $I_z$. The hyperbolic length of $\Gamma_w$ is
\[
\int_{I_z}\lambda_\Omega(x)\,dx\le\frac{3}{d_z}\cdot\frac{2\kappa d_z}{3}=2\kappa\le2,
\]
so $\Gamma_w\subset\overline{D_\D(w,2)}$, and \eqref{eq:harnack} gives $1-|u|^2\ge e^{-2}(1-|w|^2)$ for $u\in\Gamma_w$. Let $g$ be the inverse of $F|_\ell$. Then $|g'(x)|=\frac12(1-|g(x)|^2)\lambda_\Omega(x)$, and
\[
\length(\Gamma_w)=\int_{I_z}|g'(x)|\,dx\ge\frac{e^{-2}(1-|w|^2)}{2}\cdot\frac{3c_0}{4d_z}\cdot\frac{2\kappa d_z}{3}=\frac{c_0\kappa}{4e^2}(1-|w|^2).
\]
The arcs $\Gamma_w$, $w\in W$, are pairwise disjoint subsets of $F^{-1}(\R)$, by Lemma~\ref{lem:reflection}(c) and the disjointness of the intervals $I_z$. As in the proof of Lemma~\ref{lem:injectivity}, $\overline{D_\D(w,2)}\subset D(w,7(1-|w|))$. Let $\zeta\in\T$ and $r>0$. If $w\in D(\zeta,r)$, then $1-|w|<r$ and $\Gamma_w\subset D(\zeta,8r)$. Therefore
\begin{align*}
\mu_W\bigl(D(\zeta,r)\bigr)&\le\frac{4e^2}{c_0\kappa}\sum_{w\in W\cap D(\zeta,r)}\length(\Gamma_w)\le\frac{4e^2}{c_0\kappa}\length\bigl(F^{-1}(\R)\cap D(\zeta,8r)\bigr)\\
&\le\frac{32e^2}{c_0\kappa}\,r\,\|ds|_{F^{-1}(\R)}\|_{\CM}.
\end{align*}
This is the first inequality in (b), and the second follows from Corollary~\ref{cor:carleson-real} and Lemma~\ref{lem:bmo-basic}(b). The final assertion is Carleson's interpolation theorem.
\end{proof}

\subsection{Relation to the Carleson homogeneous condition}
As in the proof of Proposition~\ref{prop:interpolation}, we write $d_z=\dist(z,E)$ for $z\in\R\setminus E$. Proposition~\ref{prop:interpolation} has a converse when the sequence $Z$ also covers $\R\setminus E$ in the sense that, for some constant $0<\kappa'<1$,
\begin{equation}\label{eq:cover}
\R\setminus E\subset\bigcup_{z\in Z}(z-\kappa'd_z,\,z+\kappa'd_z).
\end{equation}

\begin{proposition}\label{prop:interp-converse}
Let $0<\kappa'<1$, let $Z\subset\R\setminus E$ be a countable set satisfying \eqref{eq:cover}, and set $R=-2\log(1-\kappa')$, $W=F^{-1}(Z)$ and $\mu_W=\sum_{w\in W}(1-|w|^2)\delta_w$. Then
\[
\|ds|_{F^{-1}(\R)}\|_{\CM}\le Re^{2R}\|\mu_W\|_{\CM}.
\]
\end{proposition}

\begin{proof}
We may assume that $\|\mu_W\|_{\CM}<\infty$. Let $x\in\R\setminus E$, and choose $z\in Z$ with $|x-z|<\kappa'd_z$. The segment $[x,z]$ lies in the component of $\R\setminus E$ containing $z$, and $\dist(t,E)\ge d_z-|t-z|$ for $t\in[x,z]$. The inequality $\lambda_\Omega\le2/\dist(\cdot,E)$ holds for every domain (see the proof of Lemma~\ref{lem:density}), and therefore
\[
\int_{[x,z]}\lambda_\Omega(t)\,dt<\int_0^{\kappa'd_z}\frac{2\,ds}{d_z-s}=R .
\]
Now let $y\in F^{-1}(\R)$, set $x=F(y)$, choose $z$ as above, and let $\ell$ be the component of $F^{-1}(\R)$ through $y$. By Lemma~\ref{lem:reflection}(c), $\ell$ contains a point $w\in W$ with $F(w)=z$, and $F$ maps the arc of $\ell$ between $y$ and $w$ onto $[x,z]$. Since $\ell$ is a geodesic and $F$ is a local isometry, $d_\D(y,w)<R$. Consequently, $F^{-1}(\R)$ is covered by the geodesic arcs $\Gamma_w=\ell_w\cap D_\D(w,R)$, $w\in W$, where $\ell_w$ is the component of $F^{-1}(\R)$ through $w$. The hyperbolic length of $\Gamma_w$ is at most $2R$, and $1-|u|^2\le e^R(1-|w|^2)$ for $u\in \Gamma_w$ by \eqref{eq:harnack}. Hence $\length(\Gamma_w)\le Re^R(1-|w|^2)$.

Let $\zeta\in\T$ and $r>0$. If $u\in \Gamma_w\cap D(\zeta,r)$, then $1-|u|<r$ and $w\in D_\D(u,R)=\sigma_u(D(0,\tanh(R/2)))$, so that
\[
|w-u|\le\frac{2(1-|u|)\tanh(R/2)}{1-\tanh(R/2)}=(e^R-1)(1-|u|)<(e^R-1)\,r .
\]
Thus $w\in D(\zeta,e^Rr)$, and
\begin{align*}
\length\bigl(F^{-1}(\R)\cap D(\zeta,r)\bigr)&\le\sum_{w\in W\cap D(\zeta,e^Rr)}\length(\Gamma_w)\le Re^R\mu_W\bigl(D(\zeta,e^Rr)\bigr)\\
&\le Re^{2R}\|\mu_W\|_{\CM}\,r .\qedhere
\end{align*}
\end{proof}

The qualitative connection between the Hayman--Wu property and interpolation of lifted Whitney sequences is already present in Fern\'andez \cite[Section 4]{Fe}. Here we give direct quantitative estimates under the separation and covering conditions \eqref{eq:sparse} and \eqref{eq:cover}.

\begin{theorem}\label{thm:interp-char}
Let $0<\kappa\le1$ and $0<\kappa'<1$, and let $Z\subset\R\setminus E$ satisfy \eqref{eq:sparse} and \eqref{eq:cover}. Then the following conditions are equivalent.
\begin{enumerate}[label=\textup{(\alph*)}]
\item $\Omega$ is a Carleson-homogeneous Denjoy domain.
\item $F^{-1}(Z)$ is an interpolating sequence for $H^\infty$.
\item $\sum_{w\in F^{-1}(Z)}(1-|w|^2)\delta_w$ is a Carleson measure on $\D$.
\end{enumerate}
\end{theorem}

\begin{proof}
If (a) holds, then $\log F'\in\BMOA(\D)$ by \cite[Theorem~1.1]{HZ}, and Proposition~\ref{prop:interpolation} gives (b). The implication (b)$\Rightarrow$(c) is part of Carleson's interpolation theorem. If (c) holds, Proposition~\ref{prop:interp-converse} shows that arclength on $F^{-1}(\R)$ is a Carleson measure on $\D$, and (a) follows from the equivalence of (i) and (iii) in Theorem~\ref{thm:main}, proved in Section~\ref{subsec:HW}.
\end{proof}

Sequences as in Theorem~\ref{thm:interp-char} always exist. For each component $(a,b)$ of $\R\setminus E$, set $h=(b-a)/2$ and take the two sequences
\[
a+2^{-k}h,\qquad b-2^{-k}h,\qquad k=0,1,2,\ldots .
\]
Let $Z$ be their union over all components, counting each midpoint once. The dyadic spacing gives \eqref{eq:sparse} with $\kappa=1/2$ within each interval, and points in different intervals satisfy $|z-z'|\ge d_z+d_{z'}$. Moreover, every point of an interval is within distance $d_z/2$ of a chosen point $z$, so \eqref{eq:cover} holds with $\kappa'=3/4$. Sequences of this type occur in Fern\'andez's proof \cite[Section~4]{Fe}, and Proposition~\ref{prop:interpolation} gives the interpolation property of their lifts directly.

Some restriction on $\kappa'$ is essential, as the following example with $\kappa'=6$ shows.

\begin{example}\label{ex:fiber}
Let $K$ be the middle-third Cantor set, translated and scaled so that $K\subset[-1/2,1/2]$ and $(-1/6,1/6)$ is the largest component of $(-1/2,1/2)\setminus K$, and let $E=[-\infty,-1]\cup K\cup[1,+\infty]$. The set $Z=\{0\}$ satisfies \eqref{eq:sparse} for every $\kappa$, and since $d_0=1/6$, the interval $(-6d_0,6d_0)=(-1,1)$ contains $\R\setminus E$,
namely, $Z$ satisfies \eqref{eq:cover} for $\kappa'=6$. The sequence $W=F^{-1}(0)$ is an interpolating sequence for $H^\infty$. However, $\Omega$ is not a Carleson-homogeneous Denjoy domain.
\end{example}

\begin{proof}
The set $E$ is not Carleson-homogeneous because $|K|=0$. Let $\psi$ be a conformal map of the simply connected domain $\Cbar\setminus([-\infty,-1]\cup[1,+\infty])$ onto $\D$ with $\psi(0)=0$. This domain contains $\Omega$, so the function $\psi\circ F$ is bounded and analytic in $\D$, and its zeros are the points of $W$, all of them simple. Hence $W$ satisfies the Blaschke condition $\sum_{w\in W}(1-|w|)<\infty$. The group $G$ acts transitively on $W$ and preserves the pseudo-hyperbolic distance $\rho$, and so the product $\prod_{v\in W\setminus\{w\}}\rho(w,v)$ does not depend on $w\in W$. This product is positive, because none of its factors vanishes and
\[
\sum_{v \in W\setminus \{w\}}\big(1-\rho(w,v)^2\big)
=\sum_{v \in W\setminus \{w\}}\frac{(1-|w|^2)(1-|v|^2)}{|1-\overline wv|^2}
\le\sum_{v \in W}\frac{1+|w|}{1-|w|}\,(1-|v|^2)<\infty.
\]
Thus $W$ is uniformly separated, and Carleson's interpolation theorem \cite[Chapter~VII]{Ga} shows that $W$ is an interpolating sequence.
\end{proof}


\begin{thebibliography}{99}

\bibitem{AZ}K. Astala and M. Zinsmeister, {\it Teichm\"uller spaces and $\mathrm{BMOA}$}, Math. Ann. {\bf 289} (1991), 613--625.

\bibitem{BP}A.F. Beardon and Ch. Pommerenke, {\it The Poincar\'e metric of plane domains}, J. London Math. Soc. (2) {\bf 18} (1978), 475--483.

\bibitem{BJ}C.J. Bishop and P.W. Jones, {\it Harmonic measure, $L^2$ estimates and the Schwarzian derivative}, J. Anal. Math. {\bf 62} (1994), 77--113.

\bibitem{Ca}L. Carleson, {\it On $H^\infty$ in multiply connected domains}, in: Conference on Harmonic Analysis in Honor of Antoni Zygmund, Vol. II, Wadsworth, Belmont, CA, 1983, 349--372.

\bibitem{Du}P.L. Duren, {\it Theory of $H^p$ Spaces}, Pure and Applied Mathematics {\bf 38}, Academic Press, New York--London, 1970.

\bibitem{Fa}K.J. Falconer, {\it The Geometry of Fractal Sets}, Cambridge Tracts in Mathematics {\bf 85}, Cambridge University Press, Cambridge, 1985.

\bibitem{Fe}J.L. Fern\'andez, {\it Domains with strong barrier}, Rev. Mat. Iberoam. {\bf 5} (1989), 47--65.

\bibitem{Ga}J.B. Garnett, {\it Bounded Analytic Functions}, revised first edition, Graduate Texts in Mathematics {\bf 236}, Springer, New York, 2007.

\bibitem{GNR}J. Gr\"ohn, A. Nicolau and J. R\"atty\"a, {\it Mean growth and geometric zero distribution of solutions of linear differential equations}, J. Anal. Math. {\bf 134} (2018), 747--768.

\bibitem{HW}W.K. Hayman and J.-M.G. Wu, {\it Level sets of univalent functions}, Comment. Math. Helv. {\bf 56} (1981), 366--403.

\bibitem{HZ}S. Huo and M. Zinsmeister, {\it Denjoy domains and $\mathrm{BMOA}$}, J. Geom. Anal. {\bf 36} (2026), Paper No. 255, 8 pp.

\bibitem{JM}P.W. Jones and D.E. Marshall, {\it Critical points of Green's function, harmonic measure, and the corona problem}, Ark. Mat. {\bf 23} (1985), 281--314.

\bibitem{Ne}Z. Nehari, {\it The Schwarzian derivative and schlicht functions}, Bull. Amer. Math. Soc. {\bf 55} (1949), 545--551.

\bibitem{Po75}Ch. Pommerenke, {\it Univalent Functions}, Vandenhoeck \& Ruprecht, G\"ottingen, 1975.

\bibitem{Po78}Ch. Pommerenke, {\it On univalent functions, Bloch functions and $\mathrm{VMOA}$}, Math. Ann. {\bf 236} (1978), 199--208.

\bibitem{Po84}Ch. Pommerenke, {\it On uniformly perfect sets and Fuchsian groups}, Analysis {\bf 4} (1984), 299--321.

\bibitem{Z}M. Zinsmeister, {\it Espaces de Hardy et domaines de Denjoy}, Ark. Mat. {\bf 27} (1989), 363--378.

\end{thebibliography}
\end{document}